\documentclass[10pt]{article}

\usepackage[dvipdf]{graphicx}
\usepackage[latin1]{inputenc}
\usepackage{latexsym}
\usepackage{amsmath,amssymb,amsfonts,amsthm}
\usepackage{xcolor}
\usepackage{mathrsfs}
\usepackage{bbm}
\usepackage{bm}
\usepackage{enumitem}

\usepackage[numbers,comma,sort]{natbib}

\definecolor{db}{RGB}{0, 0, 130}
\definecolor{wildstrawberry}{rgb}{1.0, 0.26, 0.64}

\usepackage[colorlinks=true,citecolor=db,linkcolor=black,urlcolor=blue,pdfstartview=FitH]{hyperref}

\usepackage{pdfsync}

\definecolor{rp}{rgb}{0.25, 0, 0.75}
\definecolor{dg}{rgb}{0, 0.5, 0}
\newcommand{\R}{\mathbb{R}}
\newcommand{\T}{\mathbb{T}}
\newcommand{\D}{\mathcal{D}}

\newcommand{\Z}{\mathbb{Z}}
\newcommand{\N}{\mathbb{N}}
\newcommand{\EE}{\mathbb{E}}

\newcommand{\PP}{\mathbb{P}}

\newcommand{\Ccal}{\mathcal{C}}

\newcommand{\Hyp}{\mathbf{A}}

\makeatletter
\newcommand{\customlabel}[2]{%
   \protected@write \@auxout {}{\string\newlabel {#1}{{#2}{\thepage}{#2}{#1}{}}}%
   \hypertarget{#1}{#2\hspace{-0.14cm}}
}
\makeatother

\numberwithin{equation}{section}

\newtheorem{theorem}{Theorem}[section]

\newtheorem{definition}[theorem]{Definition}%

\newtheorem{corollary}[theorem]{Corollary}

\newtheorem{lemma}[theorem]{Lemma}

\newtheorem{proposition}[theorem]{Proposition}
\newtheorem{remark}[theorem]{Remark}

\author{Christian Olivera\footnote{Departamento de Matem\'atica, Universidade Estadual de Campinas, Brazil. \texttt{colivera@ime.unicamp.br}.} \and 
Alexandre Richard\footnote{Universit\'e Paris-Saclay, CentraleSup\'elec and CNRS FR-3487. \texttt{alexandre.richard@centralesupelec.fr}.}
\and Milica Toma\v sevi\'c\footnote{CMAP, Ecole polytechnique, CNRS, I.P. Paris, 91128 Palaiseau, France. ~\hfill ~\newline\texttt{milica.tomasevic@polytechnique.edu}.
}}

\title{ \Large{\textbf{Quantitative approximation of the Burgers and Keller-Segel equations
by moderately interacting particles}}}
\begin{document}

\maketitle

\begin{abstract}
In this work we obtain rates of convergence for two moderately interacting stochastic particle systems with singular kernels associated to the viscous Burgers and Keller-Segel equations. \\
The main novelty of this work is to consider a non-locally integrable kernel. Namely for the viscous Burgers equation in $\R$, we obtain almost sure convergence of the mollified empirical measure to the solution of the PDE in some Bessel space with a rate of convergence of order $N^{-1/6}$, on any time interval. With the same rate, convergence also holds for the genuine empirical measure in Wasserstein distance, and at the level of the trajectories of the particles with the standard coupling to McKean-Vlasov particles. 

In the case of the Keller-Segel equation on a  $d$-dimensional torus, we obtain almost sure convergence of the mollified empirical measure to the solution of the PDE in some $L^q$ space with a rate of order $N^{-\frac{1}{2(d+1)}}$. The result holds up to the maximal existence time of the PDE, for any value of the chemo-attractant sensitivity $\chi$. 
\end{abstract}

\section{Introduction}

In this work we consider the following nonlinear Fokker-Planck equation:
\begin{equation}\label{eq:PDE}
\begin{cases}
&\partial_t u(t,x) = \Delta u(t,x) -  \nabla\cdot \left(u(t,x) (K\ast u )(t,x)\right),\quad t>0,~x\in\mathcal{D},\\
& u(0,x) = u_0(x), \quad x\in\mathcal{D},
\end{cases}
\end{equation}
where $\mathcal{D}$ is either $\R$ or the $d$-dimensional torus $\T^d = \R^d/\Z^d$. We focus on two particular cases with singular interaction kernel $K$. Firstly, we consider the viscous Burgers equation on $\R$ with $K = \frac{1}{2}\delta_0$, where $\delta_0$ denotes the Dirac distribution. Secondly, we treat the parabolic-elliptic Keller-Segel equation on the torus in dimension $d\geq 2$ for which $ K(x)= -\chi_d \frac{x}{ |x|^d}$, with $\chi_d>0$. 
Motivated by numerical simulation, our main objective is the quantitative approximation of \eqref{eq:PDE} by means of a stochastic particle systems in moderate interaction. 

~
 
For both equations, we study the following moderately interacting system of stochastic particles (in the sense of Oelschl\"ager~\cite{Oelschlager85} and M\'el\'eard and Roelly~\cite{Meleard}): 
\begin{equation}\label{eq:IPS0}
dX_{t}^{i,N}= \frac{1}{N}\sum_{k=1}^{N} (K \ast V^{N})(X_{t}^{i,N} -X_{t}^{k,N})\; dt + \sqrt{2} \; dW_{t}^{i}, \quad t\leq T,\quad  1\leq i \leq N,
\end{equation}
where $(W^i)_{1\leq i\leq N}$ are independent standard Brownian motions and $V^{N}(x):=N^{d\alpha}V(N^{\alpha}x), $ for some $ \alpha \in (0,1]$ and some smooth probability density $V$. For this particle system, our objective is to quantify the convergence of the empirical measure $\mu^N= (\mu^N_t)_{t\leq T}$ and the regularised empirical measure $u^N:= V^N\ast \mu^N$ towards the solution of \eqref{eq:PDE}. In the Burgers case, mollifying the interaction seems to be the only way to define the particle systems. In the Keller-Segel case, while it is possible to study the non-mollified particle system for some values of $\chi_{d}$ (\cite{FournierJourdain,CattiauxPedeches,FournierTardy,Tardy}), the mollification avoids discussions about the well-posedness of the particles and their collisions to study directly their convergence.

First, we obtain the global well-posedness for mild solutions of the viscous Burgers equation, i.e. for the PDE \eqref{eq:PDE} with $K=\frac{1}{2}\delta_{0}$, in Bessel spaces $H^\beta_{p}(\R)$ when $\beta-1/p>0$ (see Theorem~\ref{th:PDEB}). Then our  main result is that for any $T>0$,
\begin{equation}\label{eq:convergenceBurgersIntro}
\displaystyle \lim_{N\to +\infty} \, N^{\xi-\varepsilon} \sup_{t\in [0,T]} \| u^N_t-u_{t} \|_{\beta,p} =0 ~a.s.,
\end{equation}
where the rate $\xi-\varepsilon$ is in the best case equal to $\frac{1}{6}-\varepsilon$ (see Theorem~\ref{th:burgers} and Remark~\ref{rk:2}). In addition, we prove that the same rate holds for the non-mollified empirical measure $\mu^N$ in Wasserstein distance. We comment later on the optimality of this rate of convergence, see the end of this paragraph and Remark~\ref{rk:optimality}.\\
The approach to prove this result is to write the mild equation satisfied by $u^N$ and then write the difference $u-u^N$, for the mild solution $u$ of \eqref{eq:PDE}. This difference involves several terms, including a nonlinear quadratic term and a stochastic convolution integral, which create two main difficulties:
\begin{itemize}
\item First, the nonlinearity prevents us  from quantifying the convergence \emph{via} a direct Gr\"onwall argument. Instead, we propose here to use Bihari's inequality. However, this approach works only on a small random time interval. Thus a delicate point is to achieve a global-in-time  convergence, by controlling how this random time horizon behaves with $N$. In addition, we need H\"older continuity for $u^N$ in space, uniformly in $N$. Obtaining estimates in $H^\beta_{p}(\R)$ with $\beta-1/p>0$ ensures this property by Sobolev embedding.

\item The second difficulty is related to the stochastic convolution integral, see \eqref{eq:defMB}. Due to the mild formulation of the problem, this stochastic integral is function-valued and we need to bound its moments in $H^\beta_{p}(\R)$. To do so, we use Proposition~\ref{prop:martingale-bound}, which itself relies on a BDG inequality for infinite-dimensional martingales~\cite{vNeervenEtAl}.

\end{itemize}
The rate of convergence $\xi$ in \eqref{eq:convergenceBurgersIntro} appears as a trade-off between two competing quantities: on the one hand, an error appears due to the convolution with $V^N$, in the form $ \lVert \langle V^N(x-\cdot) \mu, v(x) - v(\cdot)\rangle \rVert_{L^p(\mu)}$ with $\mu$ a probability measure and $v\in H^\beta_{p}(\R)$;
and on the other hand, it accounts for the rate at which the stochastic integral that appears in $u^N$ vanishes. We obtain optimal rates of convergence for these two quantities, but we do not know if the rate from \eqref{eq:convergenceBurgersIntro} is itself optimal. We comment further on this in Remark~\ref{rk:optimality}.\\
Combining the result \eqref{eq:convergenceBurgersIntro} and some regularity obtained on Burgers PDE (see Appendix~\ref{app:regPDE}), we obtain a similar quantitative result on the trajectories of the particles (Proposition~\ref{prop:BurgersPoC-particles}), which in turn implies the trajectorial propagation of chaos, that is, convergence of the empirical measure on the space of trajectories (Corollary~\ref{cor:trajectorialPoC-Burgers}).

Secondly, we adapt the above methodology to the Keller-Segel equation and obtain an almost sure rate of convergence for the particle system. Namely, 
for the kernel $K(x)= -\chi_d \frac{x}{ |x|^d}$ on the $d$-dimensional torus (for any $d\geq 2$), we obtain (see Theorem~\ref{TMain}) that 
\begin{equation*}
\lim_{N\to +\infty} \, N^{\varrho-\varepsilon} \sup_{t\in [0,T]} \| u^N_t-u_{t} \|_{L^{q}(\mathbb{T}^d)} = 0~~ a.s.,
\end{equation*}
where $q>d$, $T$ is smaller than the maximal existence time for the PDE (see Definition \ref{def:defMild}) and $\varrho-\varepsilon$ is in the best case $\frac{1}{2(d +1)}-\varepsilon$ (see Remark~\ref{rk:1}). This rate is the same as the one found in \cite{ACM}, however no cut-off needs to be applied to the drift in the present work to obtain a.s. convergence.
This also takes into account that the rate of convergence of the initial condition, $\| u^N_0-u_{0} \|_{L^{q}(\mathbb{T}^d)}$ is at least as fast as $\frac{1}{2(d +1)}$; we give conditions to achieve such rate in Proposition~\ref{prop:initialrate}. This general result could be of independent interest, as it quantifies the approximation of probability density functions by i.i.d. samples of random variables, in Bessel norm.
Note that we presented the convergence for the Keller-Segel model in $L^q(\T^d)$ space, but convergence in Bessel spaces, like for the Burgers equation, is also possible (see Remark~\ref{rk:BesselforKS}) with the necessary adjustment on the rate of convergence.\\
Moreover, as for the Burgers equation, trajectorial convergence results are deduced (see Section~\ref{sec:PoC-KS}).

The difference between the two functional frameworks, $H^\beta_p(\R)$ for Burgers and $L^q(\T^d)$ with $q>d$ for Keller-Segel, can be roughly explained as follows.
When $K$ is the Keller-Segel kernel, it regularises enough to ensure that $K\ast u^N$ is H\"older continuous in $L^q$, with a H\"older coefficient that determines the rate of convergence $\varrho$. This is no longer the case when $K$ is a Dirac. Nevertheless, assuming initial conditions in $H^\beta_{p}(\R)$, we are able to prove that $u^N$ is uniformly bounded in $H^\beta_{p}(\R)$ and is therefore $(\beta-1/p)$-H\"older continuous. In both cases, this uniform H\"older continuity is used in our Gr\"onwall/Bihari argument to obtain a rate of convergence.

These two functional frameworks are well suited for numerical approximation. In particular, as $H^\beta_{p}(\R)$ embeds into $L^\infty$ for $\beta-1/p>0$, this could lead to pointwise numerical errors for approximation of the PDE by the empirical measure of stochastic particle systems. The work \cite{Cazacu} addresses this question, providing a rate of convergence for the Euler scheme of Coulomb-type particle systems towards the solution of the Fokker-Planck equation.

It is also worth noticing that, if one were to repeat our proofs for a Lipschitz continuous and bounded kernel, the rate of convergence would be $N^{-\frac{1}{2}-\varepsilon}$ for $d=1$ and $N^{-\frac{1}{d}-\varepsilon}$ for $d\geq 2$ in $L^1$-norm. We see that compared to the regular case, the two kernels we consider here slow down the speed, due to their singularity and the additional functional framework one needs to consider in order to tame the singularity.

\paragraph{Literature.} The following discussion is structured into four paragraphs: the first one focuses on general literature about moderately interacting particle systems; the second and third focus on quantitative convergence for respectively the viscous Burgers equation and the Keller-Segel equation, no matter whether the regime is moderate or not. The last paragraph emphasises the main differences with our previous work \cite{ACM}.

Moderately interacting particle systems with regular coefficients and their trajectorial propagation of chaos were studied initially in \cite{Oelschlager85,Meleard}.  Furthermore, rates of convergence and some types of fluctuations theorems for moderate particle systems where obtained  by \citet{Oelschlager87} and more recently \citet{HolzingerEtAl} for some specific singular kernels, or for more general regular kernels by \citet{JourdainMeleard}.
Based on a mild formulation of the empirical measure of a moderately interacting system and semigroup theory, Flandoli, Leimbach and Olivera \cite{FlandoliLeimbachOlivera} recently developed a technique to approximate nonlinear PDEs by smoothed empirical measures in strong functional topologies. This technique was also applied for a PDE-ODE system related to aggregation phenomena, see Flandoli and Leocata \cite{FlandoliLeocata}; for non-local conservation laws, see Olivera and Simon \cite{Simon}; for the 2d Navier-Stokes equation, see Flandoli, Olivera and Simon \cite{FlandoliOliveraSimon}, etc. In \cite{ACM} we further obtain convergence rates and propagation of chaos for singular repulsive and attractive Riesz potentials (including the Coulomb case).

In most of the previous models (at the exception of the porous medium equation \cite{OelschlagerPorous,HolzingerEtAl-Porous}), the interaction kernel $K$ always had to be at least locally integrable, hence excluding the Burgers interactions ($K = \frac{1}{2}\delta_0$). 
However, particle approximation of the Burgers equation has a rich history and propagation of chaos were
considered first by \citet{McKean}, then \citet{Calde}, \citet{Kac2}, \citet{Oelschlager85}, \citet{Osada}, \citet{SzTi} and \citet{BossyTalay0,BossyTalay}. However, the problem of quantifying the distance between the particles and the equation was only considered in \cite{BossyTalay0,BossyTalay}. One way to handle $K = \frac{1}{2}\delta_0$ is to 
 take as initial condition a cumulative distribution function (c.d.f.). Indeed,~\citet{BossyTalay0} transformed the equation into an integrated version of \eqref{eq:PDE}, where now the kernel becomes the integral of a Dirac, i.e. a Heavyside function. In such framework, they proved quantitative propagation of chaos~\cite{BossyTalay0} and a rate of convergence of the Euler scheme~\cite{BossyTalay} to the viscous Burgers PDE.
 In this paper, we  combine the viewpoint of moderately interacting particles with semigroup techniques, to obtain a quantitative convergence of the empirical measure in Bessel norm for $L^1(\R)$ initial conditions (i.e. without assuming the initial condition is a c.d.f.), which seems to be new in the literature.

In the second part of this work, we consider the Keller-Segel kernels ($ K(x)= -\chi_d \frac{x}{ |x|^d}$) as well as other kernels with similar scaling properties (see Remark~\ref{rk:Coulomb}), which are still singular, but locally integrable kernels.  
The Keller-Segel system of PDEs has been extensively studied for its property that a blow-up in finite time may occur in the equation. That is, the measure solution develops, in finite time, a singular part. For instance in $\R^2$, it is well known that if $\chi_2<8\pi$, the solution to \eqref{eq:PDE} exists globally (in time), while when $\chi_2>8\pi$ a blow-up in finite time occurs (see the survey of \citet{perthame},  or more recently Biler~\cite{Biler} and the references therein). 
 From the probabilistic side, the particle system (without smoothing) related to Keller-Segel equation in $\R^2$ was studied in \cite{FournierJourdain, CattiauxPedeches} and more recently in \cite{FournierTardy} and \cite{Tardy}. The authors prove well-posedness of the particle system and tightness-consistency for its empirical measure for all subcritical values of the parameter $\chi_{2}$ \cite{FournierJourdain, Tardy} and analyse particle collisions in supercritical case \cite{FournierTardy}. Furthermore, 
  Bresch, Jabin and Wang \citep{BJW2020} study
   quantitative convergence, when $N \to \infty$, of the Liouville's equations associated to $k$ fixed particles at a time $t$ towards $u_t^{\otimes k}$, where $u$ solves \eqref{eq:PDE}. Under the condition that $\chi_{2} < 8\pi$ and assuming that $u\in L^\infty((0,T); W^{2,\infty}(\T^2))$, they proved using new techniques of relative entropy the above convergence with a rate in $L^\infty((0,T); L^1((\T^2)^k))$. 
 In addition to the above references, quantitative convergence for various singular interactions has been a very active field lately, and we further mention \cite{Serfaty,JabinWang,ChodronEtAl,JabirEtAl,WangEtAl}.
   
   \medskip

To finish the discussion, let us note here that in \cite{ACM}, we obtained a rate of convergence for the empirical measure of the mollified particle system, for various kernels including the Keller-Segel one, but applying a cutoff on the drift term of each particle (i.e. we had $F(K\ast u^N)$ in the drift for some smooth cutoff $F$, instead of $K\ast u^N$ as in \eqref{eq:IPS0}). 
It was also possible to deduce a rate of convergence, but only in probability, for the particle system without cutoff \cite[Corollary 1.4]{ACM}. 
In the present work, the difficulty in lifting the cutoff appears at two different levels: first, when comparing the empirical measure of the particle system to the limit PDE in \cite{ACM}, one uses crucially the cut-off to somehow tame the nonlinearity when comparing the two drifts $u^N F(K\ast u^N)$ and $u K\ast u$ (see first bullet point in page 2 for the approach used in the present paper to circumvent this difficulty); second, in the estimation of the stochastic integral term, working with a cut-off is convenient as the particles then have bounded drifts, which makes some computations possible (see second bullet point in page 2).

\paragraph{Plan of the paper.} In Section~\ref{sec:prelim}, we list the notations used throughout the paper and in Subsection~\ref{subsec:prelimLemmas} we recall or prove some useful lemmas that hold in Bessel spaces both in the torus and in the Euclidean space. Section~\ref{sec:Burgers} and Section~\ref{sec:KS} follow the same organisation: in Subsection~\ref{subsec:Burgers1} (resp. Subsection~\ref{subsec:KS1}) we define rigorously the particle system, prove technical lemmas on the kernel and state the global well-posedness of the PDE for Burgers (resp. local well-posedness for Keller-Segel). Then in Subsection~\ref{subsec:Burgers2} (resp. Subsection~\ref{subsec:KS2}), we state and prove the convergence theorems. The trajectorial propagation of chaos is presented in Subsection~\ref{subsec:trajectorialPoC-Burgers} (resp. Subsection~\ref{sec:PoC-KS} for Keller-Segel).
 In Appendix~\ref{App}, we state and prove a bound on the stochastic convolution integral that appears in the expression of the mollified empirical measure of the Burgers particle system. In Appendix~\ref{app:initialLLN}, a rate of convergence in $L^p$/Bessel norm is obtained for the approximation of a probability density by i.i.d. samples, which could be of independent interest. Finally in Appendix~\ref{app:regPDE}, additional space regularity on the PDEs are derived, at the cost of a time explosion as $t\to 0$.

\section{Notations and preliminaries}\label{sec:prelim}

\subsection{Notations and definitions}

\begin{itemize}[leftmargin=*]

\item In this section, we use $\D$ to denote either $\T^d$ or $\R^d$.  In Section~\ref{sec:Burgers}, $\D$ will be  $\T^d$ and in Section~\ref{sec:KS}, $\D$ will be $\R$.

\item Whenever we consider a function on $\T^d$, we associate it with its periodic extension which is a function on $\R^d$.  

\item For $(X,d_X)$ a Polish space, we consider the space $\mathcal{C}([0,T];X)$ of continuous functions from $[0,T]$ to $X$ endowed with the distance
\begin{equation}\label{eq:defdtilde}
\widetilde{d}(f,g) = 1\wedge \sup_{t\in [0,T]} d_X(f(t) , g(t)) .
\end{equation}
Recall that this distance is topologically equivalent to the distance $d(f,g) = \sup_{t\in [0,T]} d_X(f(t) , g(t))$.

 Let us denote by $\mathcal{N}_{\delta}$ the H\"older seminorm of parameter $\delta\in(0,1]$, that is, for any function $f$ defined over $\D$:
\begin{equation}\label{eq:HolderNorm}
\mathcal{N}_{\delta}(f) := \sup_{x\neq y \in \D} \frac{|f(x) - f(y)|}{|x-y|^\delta}.
\end{equation}
The space of continuous and bounded functions on $\D$ which have finite $\mathcal{N}_{\delta}$ seminorm is the H\"older space $\Ccal^\delta(\D)$.

\item For $(X,d_X)$ a Polish space,  denote by $\mathcal{P}(X)$ the set of Borel probability measures on $X$.

Following \cite[Section 8.3]{BogachevII}, let us introduce the Kantorovich-Rubinstein metric which reads, for any two probability measures $\mu$ and $\nu$ on $\D$,
\begin{equation}\label{eq:defWasserstein}
\|\mu - \nu \|_{0} = \sup \left\{ \int_{\D} \phi \, d(\mu-\nu) \, ; ~ \phi \text{ Lipschitz  with } \|\phi\|_{L^\infty(\D)}\leq 1 \text{ and } \|\phi\|_{\text{Lip}} \leq 1 \right\} .
\end{equation}

 \item In this paper, $(e^{t\Delta })_{t\geq 0 }$ denotes the semigroup of the heat operator on $\D$. 
 That is, for $f \in {L}^p(\D)$, 
\begin{equation*}
e^{t \Delta}f (x)  
 = g_{2t}^{(\D)} \ast_{\D} f(x),
\end{equation*}
where $\ast_{\D}$ denotes the convolution on $\D$ and for any $t>0$, $g_{t}^{(\D)}$ is heat kernel given by
\begin{equation*}
g_{t}^{(\D)}(x) =
\begin{cases}
&\displaystyle\frac{1}{(2\pi t)^{d/2}} \sum_{k\in \mathbb{Z}^{d}} e^{-\frac{| x-k |^{2}}{2t}} ~\mbox{on } \T^d,\\
&\displaystyle\frac{1}{(2\pi t)^{d/2}} e^{-\frac{|x|^{2}}{2t}}~\mbox{on } \R^d.
\end{cases}
\end{equation*}
To make the notations easier, we will write $g_{t}$ for $g_{t}^{(\T^d)}$ and $g_{t}^0$ for $g_{t}^{(\R^d)}$.

\item If $u$ is a function or stochastic process defined on $[0,T]\times \D$, we will most of the time use the notation $u_{t}$ to denote the mapping  $x\mapsto u(t,x)$.

\item Depending on the context, the brackets $\langle\cdot , \cdot \rangle$ will denote either the scalar product in some $L^2$ space or the duality bracket between a measure and a function.

\end{itemize}

\paragraph{Bessel spaces.} We now introduce briefly \emph{Bessel potential spaces} on $\R^d$ and on the torus $\T^d$.

\begin{itemize}
\item For any  $\beta\in \R$ and $p\geq1$,  we denote by
$H^{\beta}_{p}(\mathbb{R}^{d})$ the space
\[H_{p}^{\beta}(\mathbb{R}^{d}):= \Big\{ u \text{ tempered distribution; }  \, \mathcal{F}^{-1}\Big( \big(1+|\cdot|^{2}\big)^{\frac{\beta}{2}}\; \mathcal{F} u(\cdot) \Big) \in  L^p(\mathbb R^d)\Big\},  \]
where $\mathcal Fu$ denotes the \emph{Fourier transform} of $u$. This space is endowed with the norm
 \begin{equation}\label{eq:Besselnorm}
 \| u \|_{\beta,p} = \Big\| \mathcal{F}^{-1}\big((1+|\cdot|^{2})^{\frac{\beta}{2}
}\; \mathcal{F} u(\cdot) \big) \Big\|_{L^p(\R^d)}. 
 \end{equation}
The space $H_{p}^\beta(\R^d)$ is 
associated to the Bessel potential operator $(I-\Delta)^\frac{\beta}{2}$ defined as (see e.g. \cite[p.180]{Triebel} for more details on this operator):
\begin{align*}%
(I-\Delta)^\frac{\beta}{2} f := \mathcal{F}^{-1}\left((1+|\cdot|^2)^{\frac{\beta}{2}} \mathcal{F}f \right) .
\end{align*}

\item Let $D(\mathbb{T}^{d})$ be the collection of all infinitely differentiable functions on 
$\mathbb{T}^{d}$. Then $D^{\prime}(\mathbb{T}^{d})$ stands for the topological dual of 
$D(\mathbb{T}^{d})$. We denote the Fourier coefficients of $u\in D^{\prime}(\mathbb{T}^{d})$ by $\hat{u}(k):= \frac{1}{(2\pi)^{d/2}} u(e^{2i\pi \langle k,\cdot\rangle})$. As there should be no risk of confusion, we will use the same notations for the Bessel norm and Bessel operator on the torus as in the whole space. Hence for any $\beta\in \mathbb R$, we define the Bessel potential operator $(I-\Delta)^{\beta/2}$ applied to $u\in D'(\T^d)$ by
\begin{equation*}
(I-\Delta)^{\frac{\beta}{2}} u(x)  =  \frac{1}{(2\pi)^{d/2}} \sum_{k\in \Z^d} (1+|k|^2)^{\frac{\beta}{2}} \, \hat{u}(k)\, e^{-2i\pi \langle k,x\rangle},
\end{equation*}
and denote
\begin{align}\label{eq:Besselnorm-torus}
\|u\|_{\beta,p} =  \left\| (I-\Delta)^{\frac{\beta}{2}} u\right\|_{L^p(\T^d)}.
\end{align} 
As in \cite[p.168]{SchmeisserTriebel}, $H_{p}^{\beta}(\mathbb{T}^{d})$ is
defined for $p\in\left(1,\infty\right)$ and $\beta\in \mathbb R$ by
\begin{align*}
H_{p}^{\beta}(\mathbb{T}^{d}):= \Big\{ u \in D^{\prime}(\mathbb T^d) ;\, \|u\|_{\beta,p}<\infty \Big\}.
\end{align*}
\end{itemize}

\subsection{Preliminary lemmas}\label{subsec:prelimLemmas}

We conclude this section by heat kernel estimates. First, we have that for any $t>0$,
\begin{align*}
\|\nabla g_{t}\|_{L^1(\T^d)} \leq \|\nabla g_{t}^0\|_{L^1(\R^d)} \leq \frac{C}{\sqrt{t}}.
\end{align*}
By a convolution inequality, we deduce that
\begin{align}\label{eq:boundHeat}
\|\nabla\cdot e^{t\Delta}\|_{L^p(\D)\to L^p(\D)} \leq \frac{C}{\sqrt{t}}.
\end{align}
\begin{lemma}
\label{lemma:heatBeta}
Let $\beta\in \R$ and $p\in (1, \infty)$. For any $t>0$,
\begin{align*}
\|g_{t}\|_{\beta,1} = \|(I-\Delta)^{\frac{\beta}{2}} g_{t}\|_{L^1(\D)} \leq C \left(1\vee t^{-\frac{\beta}{2}} \right),
\end{align*}
and
\begin{equation*}
 \|(I-\Delta)^{\frac{\beta}{2}} e^{t\Delta}\|_{L^p(\D)\to L^p(\D)} \leq C \left(1\vee t^{-\frac{\beta}{2}} \right).
\end{equation*}
\end{lemma}

\begin{proof}
These inequalities are standard when $\mathcal{D}=\R^d$, see e.g. Eq. (2.2) in \cite{ACM}. We prove them now for $\mathcal{D}=\T^d$. 

First, we notice that for $k\in \Z^d$, $\hat{g}_{t}(k) = \mathcal{F} g_{t}^0(2\pi k)$, where $\mathcal{F}$ denotes the Fourier transform of $\R^d$. Hence $\hat{g}_{t}(k) = g^0_{1}(2\pi k \sqrt{t})$ and
\begin{align*}
\|(I-\Delta)^{\frac{\beta}{2}} g_{t}\|_{L^1(\T^d)} &= \frac{1}{(2\pi)^{d/2}} \int_{\T^d} \Big| \sum_{k\in \Z^d} \left(1+|k|^2\right)^{\beta/2} g_{1}^0(2\pi k \sqrt{t}) \, e^{-2i\pi \langle k,x\rangle}  \Big| dx\\
&\leq \frac{1}{(2\pi)^d} \sum_{k\in \Z^d} \left(1+|k|^2\right)^{\beta/2} e^{-2\pi^2 |k|^2 t} .
\end{align*}
Since the mapping $t\mapsto t^{\beta/2} \sum_{k\in \Z^d} \left(1+|k|^2\right)^{\beta/2} e^{-2\pi^2 |k|^2 t} = \sum_{k\in \Z^d} \left(t+|k|^2 t\right)^{\beta/2} e^{-2\pi^2 |k|^2 t}$ is bounded on $[0,1]$ and the mapping $t\mapsto  \sum_{k\in \Z^d} \left(1+|k|^2\right)^{\beta/2} e^{-2\pi^2 |k|^2 t}$ is bounded on $(1,+\infty)$, this proves the first inequality.

Then for $f\in L^p(\T^d)$,
\begin{align*}
(I-\Delta)^{\frac{\beta}{2}} e^{t\Delta} f(x) &= \sum_{k\in \Z^d} \left(1+|k|^2\right)^\frac{\beta}{2} \hat{g}_{2t}(k)\, \hat{f}(k) \, e^{-2i\pi \langle k,x\rangle} \\
&= (2\pi)^d\, f\ast \left((I-\Delta)^{\beta/2}g_{2t}\right).
\end{align*}
Hence the second inequality follows by a convolution inequality.
\end{proof}

\begin{lemma}\label{lem:equivNorms} 
Let $\beta\in \R$ and $p\in (1, \infty)$. There exists $C>0$ such that for any distribution $f$,
\begin{equation*}
C^{-1} \| f\|_{\beta+1,p} \leq \|\nabla f\|_{\beta,p} +\| f\|_{\beta,p} \leq C\, \| f\|_{\beta+1,p} .
\end{equation*}
\end{lemma}
\begin{proof}
When $\D=\R^d$, this is the equivalence of norms described in \cite[Eq. (3) p.59]{TriebelTh} (note that the space $F^s_{p,2}$ in \cite{TriebelTh} is the Bessel space used here). 

When $\D=\T^d$, use the derivative in Fourier space and the equivalence of Bessel potential space and Triebel-Lizorkin space (defined in \cite[Theorem. (ii), p.167]{SchmeisserTriebel}), see \cite[Theorem. (v), p.168-169]{SchmeisserTriebel}.
\end{proof}

\section{The viscous Burgers equation}\label{sec:Burgers}

In this section, the domain is always $\D=\R$.

\subsection{Setting and preliminaries}\label{subsec:Burgers1}

Let us now study the viscous Burgers equation on the whole space, which corresponds to Equation~\eqref{eq:PDE} with $K= \frac{1}{2}\delta_{0}$:
\begin{equation}
\left\{ \begin{aligned} & \partial_{t}u(t,x) =\partial_{xx}^2 u(t,x)- \tfrac{1}{2} \partial_{x}(u(t,x)^{2}), \quad  x \in \R,~ t>0,
\\
&
u(0,x)=u_{0}(x),  \quad x \in \R. 
\end{aligned}  \right.  \label{Burgers}%
\end{equation}
Let $N\geq 1$. 
Let us introduce a mollifier that will be used both to regularise the interaction kernel in the  particle system and its empirical measure. Let $V:\R \to \R_+$ be a compactly supported,  smooth probability density function. For $\alpha \in [0,1]$ and any $x \in \R$, define
\begin{equation}\label{eq:defVN}
V^{N}(x):=N^{d\alpha}V(N^{\alpha}x).
\end{equation}
Below, $\alpha$ will be restricted to some interval $(0,\alpha_{0})$, see Assumption \eqref{H-B}.
 The corresponding particle system in moderate interaction reads 
\begin{equation}\label{eq:IPS-Bg}
\begin{cases}
dX_{t}^{i,N}= \displaystyle \frac{1}{2N}\sum_{k=1 }^{N}  V^{N}(X_{t}^{i,N} -X_{t}^{k,N})\; dt + \sqrt{2} \; dW_{t}^{i}, \quad t\leq T,~  1\leq i \leq N,  \\
X_{0}^{i,N},~ 1\leq i \leq N, \quad \text{are independent of } \{W^i,~1\leq i \leq N\},
\end{cases}
\end{equation}
where  $\{(W_{t}^{i})_{t\geq 0}, \; i\in\mathbb{N}\}$ is a family of independent standard $\R$-valued Brownian motions defined on a filtered probability space $\left(\Omega,\mathcal{F},(\mathcal{F}_{t})_{t\geq 0},\PP\right)$. 
Let us denote the empirical measure of $N$ particles by 
\begin{align}\label{eq:empiricalMeasure}
\mu^N_.= \frac{1}{N}\sum_{i=1}^N \delta_{X_{.}^{i,N}} ~,
\end{align}
and the mollified empirical measure by
\begin{equation}
\label{eq:mollEmpM}
u^N_{\cdot} := V^N\ast \mu^N_{\cdot} ~.
\end{equation}

\smallskip

It is well-known that when $u_{0}$ has good properties, the Cole-Hopf solution \cite{Cole} given by
\begin{equation}\label{eq:uCH}
u^{CH}(t,x) = \frac{g^0_{2t}\ast \left(u_{0} \, \exp\left(-\frac{1}{2}\int_{0}^\cdot u_{0}(z)\, dz \right) \right)(x)}{ g^0_{2t}\ast \exp\left(-\frac{1}{2}\int_{0}^\cdot u_{0}(z)\, dz \right)(x) }
\end{equation}
provides a classical solution to the Burgers equation~\eqref{Burgers}. For instance, if $u_{0}$ is a bounded probability density, \citet{SzTi} proved that the Cole-Hopf function solves \eqref{Burgers}. For our purpose, we will need some space regularity of the solutions. This will be achieved through the following mild notion of solution:
\begin{definition}\label{def:defMildBurguers}
Let $p > 1$ and  $\beta >1/p $. Given $u_0 \in  H^{\beta}_{p}(\R)$ and $T>0$,
 a function $u$ on $[0,T] \times \R$ is said to be a mild solution to \eqref{Burgers} on $[0,T]$ if
\begin{enumerate}[label=(\roman*)]
\item $u\in \Ccal([0,T]; H^{\beta}_{p}(\R)) $; 
\item $u$ satisfies the integral equation
\begin{equation}
\label{eq:mildBg}
u_{t} =  e^{t\Delta} u_0 - \frac{1}{2} \int_0^t \partial_{x} e^{(t-s)\Delta }  u_{s}^{2}\ ds, \quad 0 \leq t \leq T.
\end{equation}
\end{enumerate}
\end{definition}

We have the following well-posedness result for the Burgers equation.

\begin{theorem}\label{th:PDEB}
Let $p > 1$, $\beta \in (1/p,1)$, $T>0$ and $u_{0}$ such that $u_0 \in  L^1\cap H^{\beta}_{p}(\R)$. Then there exists a unique mild solution $u$ to \eqref{Burgers} in $\Ccal([0,T]; H^{\beta}_{p}(\R))$. Besides, $u$ coincides with the Cole-Hopf solution $u^{CH}$ and is unique in the larger space $L^\infty([0,T];L^\infty(\R))$.
\end{theorem}

\begin{proof}
We proceed in three steps. In Step $1$, we prove local well-posedness in the sense of Definition~\ref{def:defMildBurguers} using a fixed-point argument. In Steps $2$ and $3$, we prove that this solution coincides with $u^{CH}$ and that $\sup_{t\in [0,T]} \|u^{CH}_{t}\|_{\beta,p}<\infty$ for any $T>0$. Hence a mild solution exists for any time horizon $T$.

\paragraph{Step $1$.} Let $T>0$ to be fixed later. We set $B_{t}[u,v](x)= \frac{1}{2} \int_0^t \partial_{x} e^{(t-s)\Delta}  (u_{s} v_s)(x)\, ds$, for $u,v \in \Ccal([0,T]; H^{\beta}_{p}(\R))$, $t\in [0,T]$ and $x\in \R$. 
As the operator $(I-\Delta)^{\beta/2}$ commutes with $\partial_{x} e^{t\Delta}$, we have
\begin{align*}
\|B_{t}[u,v] \|_{\beta,p} &\leq \frac{1}{2} \int_{0}^t \left\| \partial_{x} e^{(t-s)\Delta}  (u_{s} v_s)\right\|_{\beta,p}\, ds\\
&\leq \frac{1}{2}\int_{0}^t \left\| \partial_{x} e^{(t-s)\Delta} (I-\Delta)^{\beta/2} (u_{s} v_s)\right\|_{L^p(\R)}\, ds.
\end{align*}
Hence in view of \eqref{eq:boundHeat}, we have
\[
\sup_{t\in [0,T]} \|B_{t}[u,v] \|_{\beta,p} \leq C \sup_{t\in [0,T]}  \int_0^t \frac{\|u_sv_s\|_{H^{\beta}_{p}(\R) }}{\sqrt{t-s}}  ds. 
\]
Now in view of \cite[Corollary 2.86]{BCD} and since $\beta-1/p>0$, $H^{\beta}_{p}(\R)$ is an algebra.
Hence we get
\[
\sup_{t\in [0,T]} \|B_{t}[u,v] \|_{\beta,p} \leq C  T^{1/2} \sup_{t\in [0,T]} \|u_{t}\|_{\beta,p} \, \sup_{t\in [0,T]}\|v_{t}\|_{\beta,p}  . 
\]
Thus by a standard contraction principle, see for instance Theorem 13.2 in \cite{Lemai}, we deduce that there exists a time $T^*>0$ and a unique mild solution $u$ to \eqref{Burgers} in $\Ccal([0,T]; H^{\beta}_{p}(\R))$ for any $T<T^*$. 

\paragraph{Step $2$.} We observe first that uniqueness of mild solutions holds in  $L^\infty([0,T]; L^\infty(\R))$ by a Gr\"onwall argument, for any $T>0$.
We know that $u^{CH}$ is a strong solution to the viscous Burgers equation, hence it is also a mild solution. In addition, it is  in $L^\infty([0,T]; L^\infty(\R))$, as we have in particular that $ \|u^{CH}_{t}\|_{L^\infty(\R)} \leq  \|u^{CH}_{0}\|_{L^\infty(\R)}$. Since there is $\Ccal([0,T]; H^{\beta}_{p}(\R)) \subset L^\infty([0,T]; L^\infty(\R))$ by a Sobolev embedding, $u$ and $u^{CH}$ both live in $L^\infty([0,T]; L^\infty(\R))$ and solve the mild equation, so they are equal. Hence it only remains to prove that a solution exists in $\Ccal([0,T]; H^{\beta}_{p}(\R))$, for any $T>0$. We now check that $u^{CH}$ provides such a solution.

If $T^*=+\infty$, then the problem is already solved. So assume that we only obtained $T^*<+\infty$ in Step 1, and let $T> T^*/2$. Then
\begin{align*}
\sup_{t\in [0,T]} \|u^{CH}_{t}\|_{\beta,p} &= \max \left( \sup_{t\in [0,T^*/2]} \|u^{CH}_{t}\|_{\beta,p}, \sup_{t\in [T^*/2,T]} \|u^{CH}_{t}\|_{\beta,p}\right) \\
&= \max \left( \sup_{t\in [0,T^*/2]} \|u_{t}\|_{\beta,p}, \sup_{t\in [T^*/2,T]} \|u^{CH}_{t}\|_{\beta,p}\right),
\end{align*}
since $u$ and $u^{CH}$ coincide on $[0,T^*/2]$ by the uniqueness proven in Step 1.
We will show in Step 3 that for any $T>T^*/2$, 
\begin{equation}\label{eq:bounduCH}
\sup_{t\in [T^*/2,T]} \|u^{CH}_{t}\|_{\beta,p}\leq \frac{C}{(T^*)^{1-1/(2p)}},
\end{equation}
with the constant $C$ that depends only on $\|u_0\|_{L^1(\R)}$. 
Thus $u^{CH} \in \Ccal([0,T]; H^{\beta}_{p}(\R))$, for any $T>0$.

\paragraph{Step $3$.}
Using the Sobolev embedding $W^{1,p}(\R) \hookrightarrow H^\beta_{p}(\R)$ (which holds as $\beta<1$), we have that for  $T'>0$,
\begin{align*}
\sup_{t\in [T^*/2,T']} \|u^{CH}_{t}\|_{\beta,p}\leq C \sup_{t\in [T^*/2,T']} \|u^{CH}_{t}\|_{W^{1,p}}.
\end{align*}
Denote $U_0(x)= e^{-\frac{1}{2}\int_0^x u_0(z) dz}$. As $u_0 \in L^1(\R)$ we have that 
\begin{equation}\label{eq:U0}
 e^{-\frac{1}{2}\|u_0\|_{L^1(\R)}} \leq U_0(x)\leq e^{\frac{1}{2}\|u_0\|_{L^1(\R)}}. 
\end{equation}
By direct computation and \eqref{eq:U0}, there is for $t\in [T^*/2,T']$
\begin{align*}
 \|\partial_x u^{CH}_{t}\|_{L^{p}(\R)} & = \left\|\frac{(\partial_x g^0_t \ast (u_0 U_0)) (g^0_t\ast U_0 )-( g^0_t \ast (u_0 U_0)) ( \partial_x g^0_t\ast U_0 )}{(g^0_t\ast U_0 )^2} \right\|_{L^{p}(\R)}\\
 &\leq   e^{\|u_0\|_{L^1(\R)}}\|(\partial_x g^0_t \ast (u_0 U_0)) (g^0_t\ast U_0 )-( g^0_t \ast (u_0 U_0)) ( \partial_x g^0_t\ast U_0 ) \|_{L^{p}(\R)}\\
 & \leq C \|U_0\|^2_\infty  \|u_0\|_{L^1(\R)} e^{\|u_0\|_{L^1(\R)}}\|\partial_x g^0_t \|_{L^{p}(\R)},
\end{align*}
from which \eqref{eq:bounduCH} holds.
\end{proof}

The restriction with respect to the parameter $\alpha$ and the hypothesis on the initial conditions of the system are given by the following assumption:

~

\noindent 
\begin{minipage}{\linewidth}
(\customlabel{H-B}{$\Hyp_{\alpha,\beta,p}$}): \hspace{0.35cm} 
Assume that $p\geq 2$ and $\beta\in (\frac{1}{p}, 1)$. 
\begin{enumerate}[label= ]
\item(\customlabel{C04-B}{$\Hyp^i_{\alpha,\beta,p}$}) The parameters $\alpha, \beta$ and  $p$,  satisfy
$$  0<\alpha<\frac{1}{2( 1+\beta -\frac{1}{p}) } .$$
\item(\customlabel{C0integ-B}{$\Hyp_{\alpha,\beta,p}^{ii}$})  Assume that $u_{0}\in L^1 \cap H^{\beta}_{p}(\R)$ with  $ \|u_{0}\|_{L^1(\R)} = 1$, and that there exists $\xi_{0}>0$ such that for any $m\geq1$,
$$
    \Big(\EE\left[  \left\|  u^N_{0} - u_{0} \right\|_{\beta,p}^m \right] \Big)^{1/m}\leq \frac{C}{N^{\xi_{0}}}.
$$
\end{enumerate}
\end{minipage}

\begin{remark}
When $X^{i,N}_{0}$ are i.i.d. with law $u_{0} \in H^{\beta'}_{p}(\R)$, $\beta'> \beta$, Proposition~\ref{prop:initialrate} gives a rate of convergence on $\big( \EE\big[ \left\|  u^N_{0} - u_{0} \right\|_{\beta,p}^m \big]\big)^{1/m}$. This rate is discussed further in Remark~\ref{rk:2}.
\end{remark}

\subsection{Convergence}\label{subsec:Burgers2}

Recall that $u$ is the global solution of \eqref{eq:mildBg} given by Theorem~\ref{th:PDEB}. In the following, $\mathcal{L}(\mu^N)$ denotes the law of $\mu^N$ on the space $\mathcal{C}([0,T]; H^{-\gamma}_{p}(\R))$, for some $\gamma>1-1/p$. Note that $\mu^N_{t}$, as a sum of Dirac measures, is in $H^{-\gamma}_{p}$ (recall that in dimension $1$, any Dirac measure is in $H_{p}^{1/p-1-\varepsilon}$). Moreover, for $m\geq 1$, $\mathcal{W}_{m}^{(\gamma,p)}$ denotes here the $m$-Wasserstein distance on the space of probability measures on $\mathcal{C}([0,T]; H^{-\gamma}_{p}(\R))$ endowed with the distance $d(f,g) = 1\wedge \sup_{t\in [0,T]} \|f_{t}-g_{t}\|_{-\gamma,p}$.

\begin{theorem}
\label{th:burgers}
Assume \eqref{H-B} holds true, and define   
\begin{equation}\label{eq:rateB}
\xi:= \min \left(\alpha(\beta-\frac{1}{p}) , \frac{1}{2} -\alpha \big(1 + \beta -\frac{1}{p} \big) , \xi_{0}\right) .
\end{equation}

Then for any $\varepsilon>0$,
\begin{enumerate}[label=(\roman*)]
    \item 
	$\displaystyle \lim_{N\to +\infty} \, N^{\xi-\varepsilon} \sup_{t\in [0,T]} \| u^N_t-u_{t} \|_{\beta,p} =0 ~a.s.$

    \item Let $\gamma>1-1/p$ and $m\geq1$.
      Then there exists $C>0$ such that for all $N\in \N^*$,
    \begin{equation*}
    \mathcal{W}_{m}^{(\gamma,p)}(\mathcal{L}(\mu^N), \delta_u) \leq C\, (N^{-(\xi-\varepsilon)} + N^{-\alpha(\gamma-\frac{1}{p'})}).
\end{equation*}
\end{enumerate}
\end{theorem}

\begin{remark}\label{rk:2}
In view of the constraint \eqref{C04-B}, this theorem gives the almost sure convergence of $u^N$ and $\mu^N$ for any $\alpha<\frac{1}{2}$ by choosing $p=2$ and $\beta$ close to $\frac{1}{2}$, provided the initial condition $u_{0}$ is in $H^\beta_{p}$ and \eqref{C0integ-B} is satisfied. On the other hand, the best possible rate of convergence one can get here is $\xi = (\frac{1}{6})^-$, by choosing $p\rightarrow +\infty$, $\beta=1^-$ and $\alpha = \frac{1}{6}$. For this to hold, one also needs $u_{0}^N$ to converge fast enough to $u_{0}$, that is with $\xi_{0}\geq \frac{1}{6}$; with $p$, $\beta$ and $\alpha$ as before, and assuming that $X^{i,N}_{0}$ are i.i.d. with law $u_{0} \in H^{\beta'}_{p}(\R)$, $\beta'\geq \beta+1$, Proposition~\ref{prop:initialrate} does ensure that $\xi_{0}\geq \frac{1}{6}$.
\end{remark}

\begin{remark}\label{rk:optimality}
The rate of convergence $\xi$ in \eqref{eq:rateB} depends on two competing quantities, as will be seen from the proof. First, the term $N^{-\alpha(\beta-1/p)}$ comes from the commutator $E_{t}$ (see proof below), which describes the rate of convergence of $\langle V^N(x-\cdot) \mu^N,  u^N(x) - u^N(\cdot)\rangle$ in $L^p$. This rate is optimal, in the sense that it quantifies the error between a quantity in $H^\beta_{p}(\R)$ and its mollification by $V^N$, namely that for any probability measure $\mu$ and any $v\in H^\beta_{p}(\R)$, $ \lVert \langle V^N(x-\cdot) \mu, v(x) - v(\cdot)\rangle \rVert_{L^p(\mu)} \lesssim N^{-\alpha(\beta-1/p)}$.\\
Secondly, the term $N^{-(\frac{1}{2} -\alpha (1 + \beta -\frac{1}{p}))}$ comes from the Bessel norm of the stochastic integral that appears in $u^N$, namely the process $M^N$ defined in \eqref{eq:defMB}. For this quantity, $N^{-(\frac{1}{2} -\alpha (1 + \beta -\frac{1}{p}))}$ is precisely what we expect, as the formal computation below shows:
\begin{align*}
\lVert M^N_{t}\rVert_{\beta,p} &\lesssim \lVert M^N_{t}\rVert_{\beta-(\frac{1}{p}-\frac{1}{2}),2} \hspace{5cm} \text{(by Sobolev embedding)}\\
&\lesssim  \Big(\frac{1}{N^2} \sum_{i=1}^N \int_{0}^t \lVert \partial_{x} e^{(t-s)\Delta} \rVert_{L^1}^2 \lVert V^N\rVert_{\beta-(\frac{1}{p}-\frac{1}{2}),2}^2 \, ds \Big)^\frac{1}{2} \quad \text{(by BDG's inequality)}\\
&\approx N^{-1/2} \, \lVert V^N\rVert_{\beta-(\frac{1}{p}-\frac{1}{2}),2} \approx N^{-(\frac{1}{2} -\alpha (1 + \beta -\frac{1}{p}))}.
\end{align*}
\end{remark}

\begin{proof}

Similarly to \cite[Eq. (2.3)]{ACM}, we obtain the following mild formulation for $x\in \R$:
\begin{equation*}
\begin{split}
u^N_t(x)=e^{t\Delta }u^N_0(x) -  \frac{1}{2} \int_0^t \partial_{x} e^{(t-s)\Delta } \langle \mu_s^N,  V^N ( x-\cdot) \,  u^N_s (\cdot) \rangle \ ds \\
- \frac{1}{\sqrt{2}N} \sum_{i=1}^N \int_0^t  e^{(t-s)\Delta} \partial_{x} V^N (x-X_s^{i,N})\, dW^i_s. 
\end{split}
\end{equation*}
Hence
\begin{align*}
u^N_{t}(x) - u_{t}(x) &= e^{t\Delta} (u^N_{0} - u_{0})(x) - \frac{1}{2}\int_{0}^t \partial_{x} e^{(t-s)\Delta } \left( \langle \mu_s^N,  V^N (x-\cdot) \,   u^N_s (\cdot) \rangle - (u_{s} (x))^2  \right)  ds\\
&\quad - \frac{1}{\sqrt{2}N} \sum_{i=1}^N \int_0^t  e^{(t-s)\Delta} \partial_{x} V^N (x-X_s^{i,N}) \, dW^i_s \\
&= e^{t\Delta} (u^N_{0} - u_{0})(x) +  \frac{1}{2} \int_0^t  \partial_{x} e^{(t-s)\Delta } \big( (u_s)^2 -  (u^{N}_s)^2\big)(x) \, ds \\
&\quad + E_{t}(x) + M^N_{t}(x), 
\end{align*}
where we have set
\begin{align}\label{eq:defMB}
E_{t}(x)&:= \frac{1}{2} \int_0^t \partial_{x} e^{(t-s)\Delta } \langle \mu^N_{s},  V^N (x-\cdot) \left( \,u^N_s(x) -  u^N_s(\cdot) \right)\rangle \, ds,\nonumber\\
M^N_{t}(x) &:=-\frac{1}{\sqrt{2}N} \sum_{i=1}^N \int_0^t  e^{(t-s)\Delta} \partial_{x} V^N (x-X_s^{i,N})\, dW^i_s. 
\end{align}
In view of the estimate \eqref{eq:boundHeat},
\begin{equation*}
\begin{split}
\| u^N_t-u_{t} \|_{\beta, p} &\leq \|e^{t\Delta }(u^N_0- u_0 )\|_{\beta, p} + C \int_0^t \frac{1}{\sqrt{t-s}} \| (u_{s})^2 - (u^{N}_{s})^2 \|_{\beta, p} \, ds\\
&\quad + \| E_{t}\|_{\beta,p}  + \|  M^{{N}}_{t}\|_{\beta, p}. 
\end{split}
\end{equation*}
Observe that for any $a,b\in \R$, one can write $a^2-b^2= 2a(a-b)-(a-b)^2$. The latter combined with the fact that $H^{\beta}_{p}(\R)$ is an algebra for $\beta>1/p$ (\cite[Corollary 2.86]{BCD}) leads to 
\begin{equation}\label{eq:decompuN-u_Burgers}
\begin{split}
\| u^N_t-u_{t} \|_{\beta,p} &\leq \|u^N_0- u_0 \|_{\beta,p} + C \int_0^t \frac{1}{\sqrt{t-s}} \| u_{s}\|_{\beta,p}   \|\, u_s - u^{N}_{s} \|_{\beta,p} \, ds\\
&\quad + C \int_0^t \frac{1}{\sqrt{t-s}}   \|\, u_s - u^{N}_{s} \|_{\beta,p}^2\, ds\\
&\quad + \| E_{t}\|_{\beta,p}  + \|  M^{{N}}_{t}\|_{\beta,p}. 
\end{split}
\end{equation}

$\bullet$ Let us focus on $E$. Using Lemma \ref{lem:equivNorms} and the positivity of $V^N$, we have 
$$ \| E_{t}\|_{\beta, p}
\leq C\int _0^t \|(I-\Delta)^{\frac{\beta+1}{2}} e^{(t-s)\Delta}\|_{L^p(\R)\to L^p(\R)} \left(\int_{\R} \langle \mu^N_{s},  V^N (x-\cdot) | u^N_s(\cdot)-   u^N_s(x)|\rangle^{p} \, dx\right)^{\frac{1}{p}} \, ds. $$
In view of Lemma \ref{lemma:heatBeta}, we get
\begin{align*}
\| E_{t}\|_{\beta, p}
&\leq C 
\int_0^{t}   \frac{1}{(t-s)^{\frac{1+\beta}{2}}} \left(\int_{\R} \langle \mu^N_{s},  V^N (x-\cdot) | u^N_s(\cdot)-   u^N_s(x)|\rangle^{p} \, dx\right)^{\frac{1}{p}} \, ds.
\end{align*}
Now, by the embedding $H^\beta_{p}(\R) \hookrightarrow \mathcal{C}^{\beta-1/p}(\R)$, we have the inequality $\left| u^N_s(\cdot)-   u^N_s(x)\right| \leq {\|u^N_s\|_{\beta,p} |\cdot-x|^{\beta-1/p}}$. Hence
\begin{align*}
\| E_{t}\|_{\beta, p} &\leq C 
\int_0^{t}   \frac{\| u_s^{N}\|_{\beta, p}}{(t-s)^{\frac{1+\beta}{2}}} \left(\int_{\R} \langle \mu^N_{s},  V^N (x-\cdot) \left|\cdot-x\right|^{\beta-1/p}\rangle^{p} \, dx\right)^{\frac{1}{p}} \, ds.
\end{align*}
Since $V$ is compactly supported, we have that $V^N (x-y) \left|y-x\right|^{\beta-1/p} \leq C N^{-\alpha(\beta-1/p)} V^N (x-y)$. This leads to 
\begin{align*}
 \|E_{t}\|_{\beta, p} &\leq \frac{C}{N^{\alpha(\beta-1/p)}} \int_{0}^t \frac{1}{(t-s)^{\frac{1+\beta}{2}}} \|u^N_{s}\|_{\beta, p}^2\, ds \\
&\leq \frac{C}{N^{\alpha(\beta-1/p)}} \int_{0}^t \frac{1}{(t-s)^{\frac{1+\beta}{2}}} \left(\|u_s\|_{\beta, p} + \|u^N_{s}-u_{s}\|_{\beta, p}\right)^2\, ds
\end{align*}
and using again the boundedness of $u$ we get
\[
\| E_{t}\|_{\beta,p} \leq \frac{C}{N^{\alpha(\beta-1/p)}} \left( 1 + \int_{0}^t \frac{1}{(t-s)^{\frac{1+\beta}{2} } } \|u^N_{s}-u_{s}\|_{\beta, p}^2\, ds \right).
\]

$\bullet$ 
Since $p\geq2$, by Sobolev embedding and Proposition \ref{prop:martingale-bound} we have for $\varepsilon$ small enough that
\begin{align}\label{eq:boundMB2}
\forall N\in \N^*,\quad \Big\| \sup_{s\in [0,t]} \| M^N_{s}\|_{\beta,p}  \Big\|_{L^m(\Omega)} 
\leq  \Big\| \sup_{s\in [0,t]} \| M^N_{s}\|_{\beta+\frac{1}{2}-\frac{1}{p},2}  \Big\|_{L^m(\Omega)} \nonumber\\
\leq C \, N^{-\frac{1}{2}(1-\alpha(1+ 2\beta+1-2/p ) + \varepsilon/2}.
\end{align}
By Borel-Cantelli's lemma, we deduce that there exists a random variable $A_{0}$ with finite moments such that almost surely,
\begin{equation*}
\sup_{s\in [0,t]} \| M^N_{s}\|_{ \beta,p}  \leq \frac{A_0}{N^{\xi-\varepsilon}} .
\end{equation*}

\smallskip

Similarly, using \eqref{C0integ-B}, there is a random variable $A_{1}$ with finite moments such that almost surely, $\|u^N_0- u_0 \|_{\beta,p} \leq A_{1} N^{-\xi_{0}+\varepsilon}$. Denote $A = \max(A_{0},A_{1},C)$, where $C$ is the constant in the upper bound on $\| E_{t}\|_{\beta,p}$.

Then gathering the previous bounds and using the fact that $\| u_{s}\|_{\beta,p}$ is bounded (Theorem~\ref{th:PDEB}), we have
\begin{equation*}
\begin{split}
\| u^N_t-u_{t} \|_{\beta,p} &\leq \frac{A}{N^{ \xi-\epsilon}} + C \int_0^t \frac{1}{\sqrt{t-s}}   \| u_s - u^{N}_{s} \|_{\beta,p}\, ds + C \int_0^t \frac{1}{(t-s)^{(1+ \beta)/2}}   \| u_s - u^{N}_{s} \|_{\beta,p}^{2}\, ds. 
\end{split}
\end{equation*}
Let $r>\frac{2}{1-\beta} $. Denote $\mathcal{U}^N_t = \| u^N_t-u_{t} \|_{\beta,p}^r$ and 
\begin{equation}
\label{eq:defANBG}
A_{N}=\frac{A^{r}}{N^{r(\xi-\epsilon)}}.
\end{equation}
 By H\"older's inequality, we have
\begin{equation*}
\begin{split}
\mathcal{U}^N_t &\leq A_{N} + C  T^{\frac{r}{2}-1} \int_0^t   \mathcal{U}^N_s \, ds + C T^{\frac{1-\beta}{2}r -1} \int_0^t  (\mathcal{U}^N_s)^{2} \, ds. 
\end{split}
\end{equation*}
Now apply Bihari's inequality (see e.g. \cite[Theorem 27]{Dragomir}) and get that
\begin{align}\label{eq:Bihari0BG}
    \mathcal{U}^N_t \leq G^{-1} \left( G(A_N) + C_T t \right), ~\text{ for any } t\in [0,\tau_N(\omega)) \cap [0,T],
\end{align}
where $C_T= C \left( T^{\frac{1-\beta}{2}r -1} \vee T^{\frac{r}{2}-1} \right)$ and
\begin{align*}
    &G(x) = \int_1^x \frac{1}{y+y^2} \, dy = \log\left(\frac{2x}{1+x}\right),\\
    & G^{-1}(x) = \frac{e^x}{2-e^x}, \\
    & \tau_N(\omega)= \frac{1}{C_T} \, \log \left( \frac{1+A_N(\omega)}{A_N(\omega)} \right).
\end{align*}
Hence, \eqref{eq:Bihari0BG} reads
\begin{align}\label{eq:BihariBG}
   \| u^N_t-u_{t} \|_{\beta,p}^r \leq  e^{C_T t} \frac{2A_N/(1+A_N)}{2 - e^{C_T  t } 2A_N/(1+A_N )}  
      , ~\text{ for any } t\in [0,\tau_N) \cap [0,T].
\end{align}
Let $\tilde{\Omega}$ be a measurable subset of $\Omega$ of measure $1$ on which $A$ is finite and \eqref{eq:Bihari0BG} holds. For $\omega\in \tilde{\Omega}$, there exists $N_{0}(\omega)$ such that for any $N\geq N_{0}(\omega)$, $\tau_{N_{0}}>T$. Thus 
\begin{align*}
\sup_{N\geq N_{0}(\omega)} N^{r(\xi-2\varepsilon)} \sup_{t\in [0,T]}  \| u^N_t-u_{t} \|_{\beta,p}^r \leq \sup_{N\geq N_{0}(\omega)} N^{r(\xi-2\varepsilon)} e^{C_T t} \frac{2A_N/(1+A_N)}{2 - e^{C_T  t } 2A_N/(1+A_N )}  <\infty.
\end{align*}
Hence 
$$\limsup_{N\to+\infty} N^{r(\xi-2\varepsilon)} \sup_{t\in [0,T]}  \| u^N_t-u_{t} \|_{\beta,p}^r \leq \limsup_{N\to+\infty} N^{r(\xi-2\varepsilon)} e^{C_T t} \frac{2A_N/(1+A_N)}{2 - e^{C_T  t } 2A_N/(1+A_N )} = 0 ,$$ 
which gives point $(i)$ of the theorem.

\paragraph{Proof of $(ii)$.} 
Recall that $\mathcal{W}_{m}^{(\gamma,p)}$ denotes here the Wasserstein distance on the space $\mathcal{C}([0,T]; H^{-\gamma}_{p}(\R))$ which is endowed with the distance $d(f,g) = 1\wedge \sup_{t\in [0,T]} \|f_{t}-g_{t}\|_{-\gamma,p}$, for some $\gamma>1-\frac{1}{p}$. 
As the proof is the same independently of the value of $m$, we fix $m=1$.
It comes by choosing the trivial coupling $\mathcal{L}(\mu^N)\otimes\delta_u$ that
\begin{align}
\label{eq:W1muNdeltaUBG}
\mathcal{W}^{(\gamma,p)}_1(\mathcal{L}(\mu^N), \delta_u)&\leq  \EE\Big[ 1\wedge 
\sup_{t\in[0,T]} \| \mu_{t}^N - u_{t} \|_{-\gamma,p} \Big] \nonumber\\
 &\leq \EE\Big[1\wedge \sup_{t\in[0,T]} \| u_{t}^N - u_{t}\|_{-\gamma,p} \Big] +   \EE\Big[1\wedge \sup_{t\in[0,T]} \| \mu_{t}^N - u^N_{t}\|_{-\gamma,p} \Big] \nonumber \\
&\leq  C \EE\Big[1\wedge \sup_{t\in[0,T]} \| u_{t}^N - u_{t}\|_{\beta,p} \Big] +  \EE\Big[1\wedge \sup_{t\in[0,T]} \sup_{ \|\phi\|_{\gamma,p'}\leq 1} \langle \mu_{t}^N - u^N_{t},\phi \rangle\Big].
\end{align}
We treat the first term on the right-hand side and will prove that 
\begin{equation}\label{eq:W1-uNBG}
    \EE\Big[1\wedge \sup_{t\in[0,T]}  \| u_{t}^N - u_{t}\|_{\beta,p} \Big] \leq C N^{-(\xi-\varepsilon)}.
\end{equation}
 Consider $N_0\equiv N_0(\omega)$ the smallest integer such that $\tau_{N_0}>T + \frac{1}{C_T} \log(2)$.
 Then we get
\begin{align*}
    \EE\Big[ 1\wedge \sup_{t\in[0,T]} \|u^N_t - u_t\|_{\beta,p} \Big] \leq  \EE\Big[ \mathbbm{1}_{\{N\geq N_0\}} \sup_{t\in[0,T]} \|u^N_t - u_t\|_{\beta,p} \Big] + \PP(N_0\geq N) .
\end{align*}
The above choice of $N_0$ induces that for $N\geq N_0$, we have that $\frac{2A_N}{1+A_N} e^{C_T T}= 2e^{C_T(T-\tau_N)} \leq 1$. Hence, in view of \eqref{eq:BihariBG}, which holds true for any $t\leq T$ on the event $\{N\geq N_0\}$, we obtain
\begin{align*}
    \EE\Big[ \mathbbm{1}_{\{N\geq N_0\}} \sup_{t\in[0,T]} \|u^N_t - u_t\|_{\beta,p} \Big] &\leq e^{\frac{1}{r} C_T T} \EE\left[\Big( \frac{2A_N/(1+A_N)}{2 - e^{C_T T} 2A_N/(1+A_N )} \Big)^{\frac{1}{r}} \right]\\
    &\leq e^{\frac{1}{r} C_T T} \EE\left[\Big( \frac{2A_N}{1+A_N} \Big)^{\frac{1}{r}} \right]\\
    &\leq C \, e^{\frac{1}{r} C_T T}\, N^{-(\xi-\varepsilon)},
\end{align*}
using the definition \eqref{eq:defANBG} of $A_{N}$ and the fact that $A$ has finite moments.

\smallskip

Now we estimate $\PP(N_0\geq N)$. By the definition of $N_0$, we have that $\tau_{N_0-1} \leq T + \frac{1}{C_T}\log(2)$. Hence in view of the definition of $\tau_{N_0-1}$ and $A_N$, we deduce that
\begin{equation*}
    N_0 \leq 1 + \left( 2 e^{C_T T}-1\right)^{\frac{1}{r(\xi-\varepsilon)}} A^{\frac{1}{\xi-\varepsilon}} .
\end{equation*}
Now we get 
\begin{align*}
    \PP \left( N_0\geq N \right) &\leq \PP\left( A \geq \frac{(N-1)^{\xi-\varepsilon}}{(2e^{C_T T}-1)^{\frac{1}{r}}} \right) \\
    & \leq C \frac{\EE A^k}{(N-1)^{k(\xi-\varepsilon})} ,
\end{align*}
by the Markov inequality, for any $k\geq 1$. Hence \eqref{eq:W1-uNBG} follows.
 
 ~
 
 For the second term in the right-hand side of \eqref{eq:W1muNdeltaUBG}, we observe that $H^{\gamma}_{p'}$ embeds into the H\"older space $\mathcal{C}^{\gamma-\frac{1}{p'}}$, where we recall that $\gamma-\frac{1}{p'}>1-\frac{1}{p}-\frac{1}{p'}=0$. Hence
\begin{align*}
|\langle \mu_{t}^{N}, \phi\rangle - \langle u^N_t, \phi \rangle| &=| \langle \mu_{t}^{N}, (\phi-\phi\ast V^{N})\rangle |\\
&\leq \Big\langle \mu_{t}^{N}, \int_{\R} V(y)~ |\phi(.)-  \phi( \frac{y}{N^{\alpha}}-.) |   dy \Big\rangle \\
&\leq \frac{C \|\phi\|_{\gamma-\frac{1}{p'}}}{N^{\alpha(\gamma-\frac{1}{p'})}}.
\end{align*}
Hence $\mathcal{W}^{(\gamma,p)}_1(\mathcal{L}(\mu^N), \delta_u) \lesssim N^{-(\xi-\varepsilon)} + N^{-\alpha(\gamma-\frac{1}{p'})}$.
\end{proof}

\subsection{Trajectorial propagation of chaos}\label{subsec:trajectorialPoC-Burgers}

Based on the previous theorem, we are able to prove uniform convergence on compact time intervals of the particles towards typical independent particles, which are solutions to the McKean-Vlasov equations:
\begin{align}\label{eq:McKVBurgers}
d\overline{X}_{t}^i = \frac{1}{2} u_{t}(\overline{X}_{t}^i)\, dt + dW^i_{t}, \quad \mathcal{L}(\overline{X}_{t}^i) = u_{t},
\end{align}
with the same Brownian motions and initial conditions $X_{0}^i$ as the particle system. Note that $u\in L^\infty([0,T];L^\infty(\R))$, thus by \cite{Veretennikov} and by uniqueness of the solution to the Burgers equation (see Theorem~\ref{th:PDEB}), the previous equation has a unique strong solution. Hence comparing the particles to the McKean-Vlasov particles yields what is sometimes referred to as strong trajectorial propagation of chaos. The following result is a quantitative expression of such property.

\begin{proposition}\label{prop:BurgersPoC-particles}
With the previous notations and assuming that \eqref{H-B} holds true, then for any $\varepsilon>0$,
\begin{equation*}
\lim_{N\to +\infty} N^{\xi-\varepsilon}  \sup_{t\in [0,T]} |X_t^{i,N}-\overline{X}_t^{i}| = 0 \ a.s.
\end{equation*}
\end{proposition}

\begin{proof}
In view of \eqref{eq:IPS-Bg} and \eqref{eq:McKVBurgers}, for any $t\leq T$,
\begin{align*}
\sup_{s\in [0,t]} |X_s^{i,N}-\overline{X}_s^{i}| &= \sup_{s\in [0,t]} \big| \int_{0}^s \big(u^N_{r}(X_r^{i,N}) - u_{r}(\overline{X}_r^{i})\big)\, dr \big|\\
&\leq  \int_{0}^t \sup_{s\in [0,r]} \|u_s^{N}-u_s\|_\infty \, dr + \int_0^t |u_r(X_r^{i,N})-u_r(\overline{X}_r^{i})| \, dr.
\end{align*}
From Proposition~\ref{prop:BurgersLipschitz}, $u_{r}\in H^{1+\delta}_{p}(\R)$ for any $\delta<1+\beta$, so by choosing $\delta>1/p$, we get by Sobolev embedding that $u_{r}$ is Lipschitz with $\lVert u_{r}\rVert_{\text{Lip}}\leq C_{T} r^{-(1+\delta-\beta)/2}$. This gives
\begin{align*}
\sup_{s\in [0,t]} |X_s^{i,N}-\overline{X}_s^{i}| 
&\leq t \sup_{s\in [0,t]} \|u_s^{N}-u_s\|_\infty + C_{T} \int_0^t \frac{1}{r^{(1+\delta-\beta)/2}} |X_r^{i,N} - \overline{X}_r^{i}| \, dr.
\end{align*}
Thus Gr\"onwall's inequality yields $\sup_{s\in [0,T]} |X_s^{i,N}-\overline{X}_s^{i}| \leq C_{T} \sup_{s\in [0,T]} \|u_s^{N}-u_s\|_\infty$. By Sobolev embedding, $\|u_s^{N}-u_s\|_\infty \lesssim \|u_s^{N}-u_s\|_{\beta,p}$, so Theorem~\ref{th:burgers}$(i)$ permits to conclude this proof.
\end{proof}

The previous proposition, combined with the dominated convergence theorem, readily implies 
 the convergence of the empirical measure on the space of trajectories.
\begin{corollary}\label{cor:trajectorialPoC-Burgers}
The sequence $(\mu^N)_{N\in \N^*}$ converges in law in $\mathcal{P}(\mathcal{C}([0,T];\R))$ towards $u$.
\end{corollary}

Note that this result is not a weaker statement of Theorem~\ref{th:burgers}$(ii)$. Although it comes with no rate, it is truly trajectorial: if one thinks about the test functions used to describe the convergence in Theorem~\ref{th:burgers}$(ii)$, these would be continuous and bounded functions on $\mathcal{C}([0,T];H^{-\gamma}_{p})$, so keeping in mind that $H^{-\gamma}_{p}(\R) \supset \mathcal{P}(\R)$, it is rather the marginals of $\mu^N$ that converge; however in Corollary~\ref{cor:trajectorialPoC-Burgers}, the test functions are defined on $\mathcal{P}(\mathcal{C}([0,T];\R))$, so really at the level of trajectories.

 \begin{remark}\label{rk:Burgers}
Although we could not obtain a rate of convergence on the space of trajectories for $\mu^N$ to $u$ (e.g. in Wasserstein distance), we mention that an alternative to Theorem~\ref{th:burgers}$(ii)$ is possible to get a rate for the Wasserstein distance directly on $\mu^N_{t}$ (i.e. not on its law).
Namely, denote by $\mathcal{W}_{m}$ the $1$-Wasserstein distance on probability measures on $\R$, where $\R$ is endowed with the equivalent (to the usual topology) distance $|\cdot|\wedge 1$. Then, we get
\begin{align}\label{eq:WassersteinBurgers}
\forall N\in \N^*, \quad \sup_{t\in [0,T]} \EE[\mathcal{W}_{m}(\mu^N_{t}, u_{t})] \leq C\, N^{-(\xi-\varepsilon)} .
\end{align}
Indeed, to show that \eqref{eq:WassersteinBurgers} holds, 
consider the empirical measure $\overline{\mu}^N_{t}$ of the first $N$ McKean-Vlasov particles from \eqref{eq:McKVBurgers}. As before, we write the proof only for $m=1$ for readability.
There is, for any $t\in (0,T]$,
\begin{align*}
\mathcal{W}_{1}(\mu_{t}^N,u_{t}) \leq \mathcal{W}_{1}(\mu_{t}^N,\overline{\mu}_{t}^N) + \mathcal{W}_{1}(\overline{\mu}_{t}^N,u_{t}).
\end{align*}
For the first term, the coupling $\pi^N_{t} = \frac{1}{N} \sum_{i=1}^N \delta_{(X_{t}^{i,N},\overline{X}_{t}^i)}$ gives 
\begin{align*}
\EE[\mathcal{W}_{1}(\mu_{t}^n,\overline{\mu}_{t}^n)] &\leq \EE \left[ \Big(\frac{1}{N} \sum_{i=1}^N |X^{i,N}_{t}- \overline{X}^i_{t}| \Big)\wedge 1 \right] \\
&\leq \EE [|X^{i,N}_{t}- \overline{X}^i_{t}|\wedge 1].
\end{align*}
As in the proof of Proposition~\ref{prop:BurgersPoC-particles}, we have $\EE [|X^{i,N}_{t}- \overline{X}^i_{t}|\wedge 1] \lesssim  \EE[ \sup_{s\leq t} \|u^N_{s} - u_{s}\|_{\beta,p}\wedge 1]$, that we bound using \eqref{eq:W1-uNBG}. Hence $\EE[\mathcal{W}_{1}(\mu_{t}^n,\overline{\mu}_{t}^n)] \lesssim N^{-(\xi-\varepsilon)}$.\\
The remaining term $\EE[\mathcal{W}_{1}(\overline{\mu}_{t}^N,u_{t})]$ is bounded using  Fournier-Guillin's theorem \cite[Theorem 1]{FournierGuillin}. To apply this theorem, we need $\int_{\R} |x|^q u_{t}(dx) = \EE[|\overline{X}^i_{t}|^q]<\infty$, which is satisfied for any $q>0$ since the drift in \eqref{eq:McKVBurgers} is bounded. So we get $\EE[\mathcal{W}_{1}(\overline{\mu}_{t}^N,u_{t})]\leq C\, N^{-1/2}$ for some universal constant $C>0$, and \eqref{eq:WassersteinBurgers} follows since $\xi \leq 1/2$.
\end{remark}

\section{The Keller-Segel equation}\label{sec:KS}

\subsection{Setting and Preliminaries}
\label{subsec:KS1}

In this whole section, the domain is $\D=\T^d$ with $d\geq 2$. 
We start with a precise definition of the particle system \eqref{eq:IPS0} in the case of the Keller-Segel model. The kernel $K:\T^d\to \R^d$ of this model is the periodisation of the function defined on $[-\frac{1}{2},\frac{1}{2})^d\setminus\{0\}$ by $K_{0}(x)= -\chi_d \frac{x}{|x|^d}$. 

First, notice that $K\in L^p(\T^d)$ for any $p\in [1, \frac{d}{d-1})$.  Hence, by the H\"older inequality, it holds
\begin{align}\label{eq:Kastf}
\|K\ast f\|_{L^\infty(\mathbb{T}^d)} \leq C \|f\|_{L^{q}(\mathbb{T}^d)},
\end{align}
for any $q>d$.

 The following property of the kernel will be frequently used:
\begin{lemma}\label{lem:boundHolder}
Let $q\in (d,+\infty)$. Then we have
\begin{align}\label{eq:HolderK}
\forall f\in L^q(\T^d), \quad \mathcal{N}_{1-\frac{d}{q}}(K\ast f) \leq C\, \|f\|_{L^q(\T^d)} .
\end{align}
\end{lemma}

 \begin{proof}
 First, we will use the fact that $H^1_{q}(\T^d)$ is continuously embedded into the H\"older space of parameter $1-d/q$, denoted by $\mathcal{C}^{1-d/q}(\T^d)$. This comes from the equality between $H^1_{q}(\T^d)$ and the Triebel-Lizorkin space $F^{1}_{q,2}(\T^d)$ (see \cite[Theorem. (v), p.168]{SchmeisserTriebel}), the equality between $\mathcal{C}^{1-d/q}(\T^d)$ and the Zygmund space $\mathscr{C}^{1-d/q}(\T^d)$ (see \cite[Theorem. (i)-(ii), p.168]{SchmeisserTriebel}), and the continuous embedding $F^{1}_{q,2}(\T^d) \subset \mathscr{C}^{1-d/q}(\T^d)$ (see \cite[Corollary. (ii), p.170]{SchmeisserTriebel}). Hence we get
 \begin{align*}
 \mathcal{N}_{1-\frac{d}{q}}(K\ast f) \leq C\, \|K\ast f\|_{1,q}.
 \end{align*}
Now by Lemma \ref{lem:equivNorms}, we have that $\|K\ast f\|_{1,q} \leq C\, (\|\nabla K\ast f\|_{0,q} +  \| K\ast f\|_{0,q}) = C (\|\nabla K\ast f\|_{L^q(\T^d)} + \| K\ast f\|_{L^q(\T^d)})$, where the equality is the Littlewood-Paley theorem. We conclude using that $K\in L^1(\T^d)$, which implies by convolution inequality that $\| K\ast f\|_{L^q(\T^d)}\leq C\, \| f\|_{L^q(\T^d)}$. As for $\|\nabla K\ast f\|_{L^q(\T^d)}$, use that $K = -\chi_{d}\nabla G + \vartheta$, where $G$ is the usual Poisson potential
and $\vartheta$ is a $\mathcal{C}^\infty$ correction to make $K$ periodic. From \cite[Theorem 9.9]{GilbargTrudinger}, 
$\nabla^2 G$ is a Calder\'on-Zygmund operator and it thus satisfies $\|\nabla^2 G\ast f\|_{L^q(\T^d)} \leq C\, \lVert f\rVert_{L^q(\T^d)}$; as for the remaining term, $\|\nabla \vartheta \ast f\|_{L^q(\T^d)} \leq C\, \lVert f \rVert_{L^q(\T^d)}$ follows by the integrability of $\nabla \vartheta$.
 \end{proof}

Let us introduce a mollifier that will be used both to regularise the interaction kernel in the  particle system and its empirical measure. Let $V:\R^d \to \R_+$ be a smooth probability density function with support in $(-\frac{1}{2},\frac{1}{2})^d$. For any $x \in [-\frac{1}{2},\frac{1}{2})^d$, define
\begin{equation*}%
V_{0}^{N}(x):=N^{d\alpha}V(N^{\alpha}x), \qquad \text{for some } \alpha \in [0,1],
\end{equation*}
and let $V^N:\T^d \to \R_+$ be the periodisation of  $V_0^N$. Below, $\alpha$ will be restricted to some interval $(0,\alpha_{0})$, see Assumption \eqref{C04}.

Let $T>0$. 
For each $N\in\mathbb{N}$, the particle system \eqref{eq:IPS0} reads more precisely:
\begin{equation*}%
\begin{cases}
dX_{t}^{i,N}= \displaystyle \frac{1}{N}\sum_{k=1 }^{N} (K \ast V^{N})(X_{t}^{i,N} -X_{t}^{k,N})\; dt + \sqrt{2} \; dW_{t}^{i}, \quad t\leq T,~  1\leq i \leq N,  \\
X_{0}^{i,N},~ 1\leq i \leq N, \quad \text{are independent of } \{W^i,~1\leq i \leq N\},
\end{cases}
\end{equation*}
where  $\{(W_{t}^{i})_{t\in[0,T]}, \; i\in\mathbb{N}\}$ is a family of independent standard $\R^d$-valued Brownian motions defined on a filtered probability space $\left(\Omega,\mathcal{F},(\mathcal{F}_{t})_{t\geq 0},\PP\right)$. The particles are periodised, hence we consider $X_{t}^{i,N}$ as a $\T^d$-valued random variable. 
The empirical measure of $N$ particles and the mollified empirical measure are denoted by $\mu^N$ and $u^N$ respectively as in \eqref{eq:empiricalMeasure} and \eqref{eq:mollEmpM}.

 ~
 
We will study the convergence of $u^N$ in $L^q(\T^d)$ norm for any $q>d$. The extension to Bessel spaces, as for the Burgers equation, is discussed in Remark~\ref{rk:BesselforKS}. The restriction with respect to the parameter $\alpha$ and the hypothesis on the initial conditions of the system are given by the following assumption:

~

\noindent 
\begin{minipage}{\linewidth}
(\customlabel{H}{$\Hyp_{\alpha,q}$}): \hspace{0.35cm} 
Let $q>d$ and assume that:
\begin{enumerate}[label= ]
\item(\customlabel{C04}{$\Hyp^i_{\alpha,q}$}) The parameters $\alpha$ and  $q>d$ satisfy
$$ 0<\alpha<\frac{1}{2 d (1 - \frac{1}{q})} .$$
\item(\customlabel{C0integ}{$\Hyp_{\alpha,q}^{ii}$}) Assume that $u_{0}\in L^1 \cap L^{q}(\mathbb{T}^d)$ with $\|u_{0}\|_{L^1(\T^d)} = 1$, and that there exists $\varrho_{0}>0$ such that for any $m\geq1$,
$$
  \Big(\EE\left[  \left\| u_0^N- u_{0} \right\|_{L^{q}(\T^d)}^m \right]\Big)^{1/m} \leq \frac{C}{N^{\varrho_{0}}}. 
$$
\end{enumerate}
\end{minipage}

\vspace{0.4cm}

~

We aim to prove the convergence of the mollified empirical measure to the PDE \eqref{eq:PDE}. As \eqref{eq:PDE} preserves the positivity and total mass $M:= \int_{\T^d} u_0 (x) \, dx$, we will assume throughout the paper that $M=1$. 

Solutions to \eqref{eq:PDE} will be understood in the following mild sense:
\begin{definition}\label{def:defMild} 
Let $q>d$. 
Given $u_0 \in L^q(\T^d)$ and $T>0$,
 a function $u$ on $[0,T] \times \T^d$ is said to be a mild solution to \eqref{eq:PDE} on $[0,T]$ if
\begin{enumerate}[label=(\roman*)]
\item $u\in \Ccal([0,T]; L^q(\T^d)) $; 
\item $u$ satisfies the integral equation
\begin{equation}
\label{eq:mildKS}
u_{t} =  e^{t\Delta} u_0 -  \int_0^t \nabla \cdot e^{(t-s)\Delta }  (u_{s}\,  K \ast u_{s})\ ds, \quad 0 \leq t \leq T.
\end{equation}
\end{enumerate}
A function $u$ on $[0,\infty) \times \T^d$ is said to be a global mild solution to \eqref{eq:PDE} if it is a mild solution to \eqref{eq:PDE} on $[0,T]$ for all $T>0$, and it is said to be local otherwise.
\end{definition}

We have the following local well-posedness for the Keller-Segel equation.

\begin{proposition}\label{prop:PDE}
Let $q>d$ such that $u_{0}\in L^q(\T^d)$. Then there exists a maximal time $T^*\in (0,+\infty]$ and a unique mild solution to the Fokker-Planck equation which is in $\mathcal{C}\left([0,T]; L^q(\T^d)\right)$, for any $T<T^*$.
\end{proposition}

\begin{proof}
Denote $\mathcal{X} = \mathcal{C}\left([0,T];L^q(\T^d)\right)$ and $ \|\cdot\|_{\mathcal{X}}$ the associated norm. Then define $B(u,v)= \int_0^t \nabla \cdot e^{(t-s)\Delta }  (u_{s} K\ast v_s)\, ds$ for $u,v \in \mathcal{X}$. In view of \eqref{eq:boundHeat}, we have 
\begin{align*}
\|B(u,v) \|_{\mathcal{X}} &\leq C  \int_0^T \frac{\|u_s K\ast v_s\|_{L^q(\mathbb{T}^d)}}{\sqrt{t-s}}\,  ds\\
&\leq C  \int_0^T \frac{\|u_s\|_{L^q(\mathbb{T}^d)} \|K\ast v_s\|_{L^\infty(\T^d)} }{\sqrt{t-s}}\,  ds .
\end{align*}
Now using \eqref{eq:Kastf}, we get
\[
\|B(u,v) \|_{\mathcal{X}} \leq C \, T^{1/2}   \|u\|_{\mathcal{X}} \|v\|_{\mathcal{X}}  . 
\]
Thus by the same contraction principle as before, see Theorem 13.2 in \cite{Lemai}, we deduce the result for $T$ small enough. 
\end{proof}

Unlike the Burgers equation, we do not expect to have a global solution for any Keller-Segel kernel $K$, depending on the dimension $d$ and the intensity $\chi_{d}$ of the interaction. Thus we keep working with local solutions in the sequel.

\subsection{Convergence}\label{subsec:KS2}

 In the following, $\mathcal{L}(\mu^N)$ denotes the law of $\mu^N$ in the space $\mathcal{C}([0,T];\mathcal{P}(\mathbb{T}^d))$, where $\mathcal{P}(\mathbb{T}^d)$ is endowed with the  Kantorovich-Rubinstein distance $ \|\cdot\|_{0}$, and for $m\geq1$, $\mathcal{W}_{m}^{(\mathcal{P})}$ denotes the $m$-Wasserstein distance on the space of probability measures on $\mathcal{C}([0,T];\mathcal{P}(\mathbb{T}^d))$ endowed with the distance $d(\mu,\nu) = 1\wedge \sup_{t\in [0,T]} \|\mu_{t}-\nu_{t}\|_{0}$.

\begin{theorem}\label{TMain}
Recall that $T^*$ is the maximal time of existence of \eqref{eq:mildKS} from Proposition~\ref{prop:PDE}, and let $T<T^*$.
Assume that \eqref{H} holds and define 
\begin{equation}\label{eq:rateKS}
\varrho= \min \left(\alpha(1-\frac{d}{q}) , \frac{1}{2}- \alpha d(1-\frac{1}{q}), \varrho_{0}\right).
\end{equation}
Then for any $\varepsilon>0$,
\begin{enumerate}[label=(\roman*)]
    \item %
	$\displaystyle \lim_{N\to +\infty} \, N^{\varrho-\varepsilon} \sup_{t\in [0,T]} \| u^N_t-u_{t} \|_{L^{q}(\mathbb{T}^d)} = 0~~ a.s.$

    \item Besides, for any $m\geq1$, there exists $C>0$ such that for all $N\in \N^*$,
    \begin{equation*}
    \mathcal{W}_{m}^{(\mathcal{P})}(\mathcal{L}(\mu^N), \delta_u) \leq C N^{-(\varrho-\varepsilon)}.
\end{equation*}
\end{enumerate}
\end{theorem}

\smallskip

\begin{remark}\label{rk:Coulomb}
 This theorem also holds if one replaces the Keller-Segel kernel by an integrable kernel that satisfies the uniform bound \eqref{eq:Kastf} and the H\"older estimate \eqref{eq:HolderK}. This is for instance the case of the periodised version of the Coulomb kernel given by $K(x) = \pm \nabla |x|^{-(d-2)} $ and of the Biot-Savart kernel $K(x) = \frac{1}{\pi} \frac{x^\perp}{|x|^2}$ in dimension $2$.
\end{remark}

\begin{remark} \label{rk:1}
 In view of the constraint \eqref{C04}, this theorem gives the almost sure convergence of $u^N$ and $\mu^N$ for any $\alpha<\frac{1}{2(d-1)}$ by choosing $q$ close to $d$, provided the initial condition has enough regularity. On the other hand, the best possible rate of convergence one can get here is $\varrho = \frac{1}{2(d+1)}$, by choosing $q=+\infty$ and $\alpha = \frac{1}{2(d+1)}$. This is comparable to Theorem 1.3 in \cite{ACM}, where a cutoff was applied on the drift of the particles, but where the particles lived in $\R^d$. In Corollary 1.4 of \cite{ACM}, convergence was also obtained for the particle system without cut-off, but only in probability. For the aforementioned parameters $q$ and $\alpha$, and with i.i.d. initial conditions $X^{i,N}_{0}$ following the law $u_{0}\in H^1_{q}(\T^d)$, Proposition~\ref{prop:initialrate} implies that $\varrho_{0}\geq \frac{1}{2(d+1)}$, so that $\varrho = \frac{1}{2(d+1)}$ is indeed achieved for such initial condition.
\end{remark}

\begin{remark}
As already mentioned in the introduction, an interesting feature of Keller-Segel equations is that they exhibit a blow-up for large enough $\chi_{d}$. In \cite{SenbaSuzuki}, a slightly modified $2$-dimensional Keller-Segel equation on the torus was studied and a condition on $\chi_{2}$ for existence of weak solutions was provided. Whether there is a blow-up or not, our approach provides an approximation of the PDE until maximal existence time.
\end{remark}

\begin{proof}
We prove $(i)$ first and then use it to prove $(ii)$.
\paragraph{Proof of $(i)$.} 
Similarly to \cite[Eq. (2.3)]{ACM}, we obtain the following mild formulation for $x\in \T^d$:
\begin{equation*}
\begin{split}
u^N_t(x)=e^{t\Delta }u^N_0(x) -  \int_0^t \nabla \cdot e^{(t-s)\Delta } \langle \mu_s^N,  V^N ( x-\cdot) \,  K\ast u^N_s (\cdot) \rangle \ ds \\
- \frac{1}{N} \sum_{i=1}^N \int_0^t  e^{(t-s)\Delta} \nabla V^N (x-X_s^{i,N})\cdot dW^i_s. 
\end{split}
\end{equation*}
Hence
\begin{align*}
u^N_{t}(x) - u_{t}(x) &= e^{t\Delta} (u^N_{0} - u_{0})(x) - \int_{0}^t \nabla \cdot e^{(t-s)\Delta } \left( \langle \mu_s^N,  V^N (x-\cdot) \,  K\ast u^N_s (\cdot) \rangle - u_{s} (x) \, K\ast u_s (x) \right) \, ds\\
&\quad - \frac{1}{N} \sum_{i=1}^N \int_0^t  e^{(t-s)\Delta} \nabla V^N (x-X_s^{i,N})\cdot dW^i_s \\
&= e^{t\Delta} (u^N_{0} - u_{0})(x) +  \int_0^t  \nabla \cdot e^{(t-s)\Delta } \left(u_{s}   \, K\ast u_s - u^{N}_{s}\, K\ast u^{N}_s\right)(x) ~ ds \\
&\quad + E_{t}(x) + M^N_{t}(x), 
\end{align*}
where we have set
\begin{align}
E_{t}(x)&:= \int_0^t \nabla \cdot e^{(t-s)\Delta } \langle \mu^N_{s},  V^N (x-\cdot) \left( \, K\ast u^N_s(x) -  K\ast u^N_s(\cdot) \right)\rangle \ ds, \label{eq:defE-KS}\\
M^N_{t}(x) &:=-\frac{1}{N} \sum_{i=1}^N \int_0^t  e^{(t-s)\Delta} \nabla V^N (x-X_s^{i,N})\cdot dW^i_s. \label{eq:defMN-KS}
\end{align}
In view of the estimate \eqref{eq:boundHeat}, one has
\begin{equation*}
\begin{split}
\| u^N_t-u_{t} \|_{L^{q}(\mathbb{T}^d)} &\leq \|e^{t\Delta }(u^N_0- u_0 )\|_{L^{q}(\mathbb{T}^d)} \\
&\quad + C \int_0^t \frac{1}{\sqrt{t-s}} \| u_{s}  \, K\ast u_s - u^{N}_{s} \, K\ast u^{N}_s \|_{L^{q}(\mathbb{T}^d)} ds\\
&\quad + \| E_{t}\|_{L^{q}(\mathbb{T}^d)}  + \|  M^{{N}}_{t}\|_{L^{q}(\mathbb{T}^d)}. 
\end{split}
\end{equation*}
It follows that 
\begin{equation*}
\begin{split}
\| u^N_t-u_{t} \|_{L^{q}(\mathbb{T}^d)} &\leq \|u^N_0- u_0 \|_{L^{q}(\mathbb{T}^d)} + C 
\int_0^t \frac{1}{\sqrt{t-s}} \| K\ast u^N_{s}\|_{L^\infty(\mathbb{T}^d)} \| u^{N}_{s}-u_{s} \|_{L^{q}(\mathbb{T}^d)} ds\\
&\quad +  C \int_{0}^t \frac{1}{\sqrt{t-s}} \|u_{s}\|_{L^q(\T^d)} \|K\ast (u_s - u^{N}_{s}) \|_{L^{\infty}(\mathbb{T}^d)} ds  \\
&\quad + \| E_t\|_{L^{q}(\mathbb{T}^d)}  + \|  M^{{N}}_{t}\|_{L^{q}(\mathbb{T}^d)} .
\end{split}
\end{equation*}
Using \eqref{eq:Kastf}, we get
\begin{equation*}
\begin{split}
\| u^N_t-u_{t} \|_{L^{q}(\mathbb{T}^d)} &\leq \|u^N_0- u_0 \|_{L^{q}(\mathbb{T}^d)} + C 
\int_0^t \frac{1}{\sqrt{t-s}} \| u^{N}_{s}-u_{s} \|_{L^{q}(\mathbb{T}^d)} \left(\| u^{N}_{s}-u_{s} \|_{L^{q}(\mathbb{T}^d)} + \| u_{s} \|_{L^{q}(\mathbb{T}^d)}\right)  ds\\
&\quad +  C \int_{0}^t \frac{1}{\sqrt{t-s}} \|u_{s}\|_{L^q(\T^d)} \| u_s - u^{N}_{s} \|_{L^{q}(\mathbb{T}^d)} ds  \\
&\quad + \| E_{t}\|_{L^{q}(\mathbb{T}^d)}  + \| M^{{N}}_{t}\|_{L^{q}(\mathbb{T}^d)} .
\end{split}
\end{equation*}
We know from Proposition \ref{prop:PDE} that $u$ is bounded in $\mathcal{C}\left([0,T]; L^q(\T^d)\right)$. Hence 
\begin{equation}\label{eq:decompKSLq}
\begin{split}
\| u^N_t-u_{t} \|_{L^{q}(\mathbb{T}^d)} &\leq \|u^N_0- u_0 \|_{L^{q}(\mathbb{T}^d)} + C 
\int_0^t \frac{1}{\sqrt{t-s}} \| u^{N}_{s}-u_{s} \|_{L^{q}(\mathbb{T}^d)}^2 \, ds\\
&\quad +  C \int_{0}^t \frac{1}{\sqrt{t-s}}  \| u_s - u^{N}_{s} \|_{L^{q}(\mathbb{T}^d)} ds  \\
&\quad + \| E_{t}\|_{ L^{q}(\mathbb{T}^d)}  + \|  M^{{N}}_{t}\|_{L^{q}(\mathbb{T}^d)} .
\end{split}
\end{equation}

~

$\bullet$ Let us bound $E$. By using the positivity of $V^N$, we get
\begin{align*}
\| E_{t}\|_{L^{q}(\T^{d}) } 
&\leq C 
\int_0^{t} \frac{1}{\sqrt{t-s}} \left(\int_{\T^{d}} \langle \mu^N_{s},  V^N (x-\cdot) \left| K\ast u^N_s(\cdot)-  K\ast u^N_s(x)\right|\rangle^{q} \, dx\right)^{\frac{1}{q}} \, ds.
\end{align*}
Using  the H\"older continuity of $K\ast u^N$ from Lemma~\ref{lem:boundHolder}, we get
\begin{align*}
\| E_{t}\|_{L^{q}(\T^{d}) } &\leq C 
\int_0^{t}   \frac{\| u_s^{N}\|_{L^q(\T^d)}}{\sqrt{t-s}} \left(\int_{\T^{d}} \langle \mu^N_{s},  V^N (x-\cdot) \left|\cdot-x\right|^{1-\frac{d}{q}}\rangle^{q} \, dx\right)^{\frac{1}{q}} \, ds.
\end{align*}
Since $V$ is compactly supported,
 we have that $V^N (x-y) \left|y-x\right|^{1-\frac{d}{q}} \leq N^{-\alpha(1-\frac{d}{q})} V^N (x-y)$. Thus, 
\begin{align*}
\sup_{s\in [0,t]} \|E_{s}\|_{L^{q}(\mathbb{T}^d)} &\leq \frac{C}{N^{\alpha(1-d/q)}} \int_{0}^t \frac{1}{\sqrt{t-s}} \|u^N_{s}\|_{L^{q}(\mathbb{T}^d)}^2\, ds \\
&\leq \frac{C}{N^{\alpha(1-d/q)}} \int_{0}^t \frac{1}{\sqrt{t-s}} \left(\|u_s\|_{L^{q}(\mathbb{T}^d)} + \|u^N_{s}-u_{s}\|_{L^{q}(\mathbb{T}^d)}\right)^2\, ds
\end{align*}
and using again the boundedness of $u$ we get
\begin{equation*}
\sup_{s\in [0,t]} \|E_{s}\|_{L^{q}(\mathbb{T}^d)} \leq \frac{C}{N^{\alpha(1-d/q)}} \left( 1 + \int_{0}^t \frac{1}{\sqrt{t-s}}  \|u^N_{s}-u_{s}\|_{L^{q}(\mathbb{T}^d)}^2\, ds \right) .
\end{equation*}

~

$\bullet$ We now focus on $M^N$. 
For $q\geq 2$, we have the embedding of $H^{d(\frac{1}{2}-\frac{1}{q})}_{2}(\T^d)$ into $L^q(\T^d)$. Hence by Proposition~\ref{prop:martingale-bound}, it comes that for any $q\geq 2$ and any $m\geq 1$,
\begin{align}\label{eq:boundMN-KS}
\forall N\in \N^*,\quad \Big\| \sup_{s\in [0,t]} \| M^N_{s}\|_{L^q(\mathbb{T}^d)} \Big\|_{L^m(\Omega)} 
&\leq C \Big\| \sup_{s\in [0,t]} \| M^N_{s}\|_{H^{d(\frac{1}{2}-\frac{1}{q})}_{2}(\T^d)} \Big\|_{L^m(\Omega)} \nonumber\\
&\leq C N^{-\frac{1}{2}(1-2\alpha d (1-\frac{1}{q}))+\varepsilon/2}.
\end{align}
By Borel-Cantelli's lemma, we deduce that there exists a random variable $A_{0}$ with finite moments such that almost surely,
\begin{equation*}
\sup_{s\in [0,t]} \| M^N_{s}\|_{L^q(\mathbb{T}^d)}  \leq \frac{A_0}{N^{\frac{1}{2}(1-2 \alpha d (1-1/q))-\varepsilon}} .
\end{equation*}

Similarly, using \eqref{C0integ}, there is a random variable $A_{1}$ with finite moments such that almost surely, $\|u^N_0- u_0 \|_{L^q(\T^d)} \leq A_{1} N^{-\varrho_{0}+\varepsilon}$. Denote $A = \max(A_{0},A_{1},C)$, where $C$ is the constant in the upper bound on $\sup_{s\in [0,t]}\| E_{s}\|_{L^q(\T^d)}$.

~

$\bullet$ 
Putting altogether the previous bounds, we have that almost surely,
\begin{align*}
\| u^N_t-u_{t} \|_{L^{q}(\mathbb{T}^d)} &\leq 
\frac{A}{N^{\rho-\varepsilon}} 
+ C \int_0^t \frac{1}{\sqrt{t-s}} \left(\| u^{N}_{s}-u_{s} \|_{L^{q}(\mathbb{T}^d)}+ \| u^{N}_{s}-u_{s} \|_{L^{q}(\mathbb{T}^d)}^2\right) \, ds\\
&\leq \frac{A}{ N^{\varrho-\varepsilon} } + C 
\int_0^t \frac{1}{\sqrt{t-s}} \left(\| u^{N}_{s}-u_{s} \|_{L^{q}(\mathbb{T}^d)}+ \| u^{N}_{s}-u_{s} \|_{L^{q}(\mathbb{T}^d)}^2\right) \, ds.
\end{align*}
Denote $\mathcal{V}^N_t = \| u^N_t-u_{t} \|_{L^{q}(\mathbb{T}^d)}^3$ and 
\begin{equation}\label{eq:defAN}
A_N = \frac{A^3}{N^{3(\varrho-\varepsilon)}}. 
\end{equation}
Then
apply H\"older's inequality with the exponents $\frac{3}{2}$ and  $3$ to obtain
\begin{equation*}
\mathcal{V}^N_t \leq A_N + C 
T^{\frac{1}{2}} \int_0^t \left(\mathcal{V}^N_s+ (\mathcal{V}^N_s)^2\right) \, ds .
\end{equation*}

Now we apply Bihari's inequality (see e.g. \cite[Theorem 27]{Dragomir}) and get that
\begin{align}\label{eq:Bihari0}
    \mathcal{V}^N_t \leq G^{-1} \left( G(A_N) + C_T t \right), ~\text{ for any } t\in [0,\tau_N(\omega)) \cap [0,T],
\end{align}
where $C_T= C\sqrt{T}$ and
\begin{align*}
    &G(x) = \int_1^x \frac{1}{y+y^2} \, dy = \log\left(\frac{2x}{1+x}\right),\\
    & G^{-1}(x) = \frac{e^x}{2-e^x}, \\
    & \tau_N(\omega)= \frac{1}{C_T} \, \log \left( \frac{1+A_N(\omega)}{A_N(\omega)} \right).
\end{align*}
Hence, \eqref{eq:Bihari0} reads
\begin{align}\label{eq:Bihari}
   \| u^N_t-u_{t} \|_{L^{q}(\mathbb{T}^d)}^3 \leq  e^{C_T t} \frac{2A_N/(1+A_N)}{2 - e^{C_T  t } 2A_N/(1+A_N )}  
      , ~\text{ for any } t\in [0,\tau_N) \cap [0,T].
\end{align}
Let $\tilde{\Omega}$ be a measurable subset of $\Omega$ of measure $1$ on which $A$ is finite and \eqref{eq:Bihari0} holds. For $\omega\in \tilde{\Omega}$, there exists $N_{0}(\omega)$ such that for any $N\geq N_{0}(\omega)$, $\tau_{N_{0}}>T$. Thus 
\begin{align*}
\sup_{N\geq N_{0}(\omega)} N^{3(\varrho-2\varepsilon)} \sup_{t\in [0,T]}  \| u^N_t-u_{t} \|_{L^{q}(\mathbb{T}^d)}^3 \leq \sup_{N\geq N_{0}(\omega)} N^{3(\varrho-2\varepsilon)} e^{C_T t} \frac{2A_N/(1+A_N)}{2 - e^{C_T  t } 2A_N/(1+A_N )}  <\infty.
\end{align*}
Hence 
$$\limsup_{N\to+\infty} N^{3(\varrho-2\varepsilon)} \sup_{t\in [0,T]}  \| u^N_t-u_{t} \|_{L^{q}(\mathbb{T}^d)}^3 \leq \limsup_{N\to+\infty} N^{3(\varrho-2\varepsilon)} e^{C_T t} \frac{2A_N/(1+A_N)}{2 - e^{C_T  t } 2A_N/(1+A_N )} = 0 ,$$ 
which gives point $(i)$ of the theorem.

\paragraph{Proof of $(ii)$.} 
As the proof is the same independently of the value of $m$, we fix $m=1$.
 Recalling the definition \eqref{eq:defWasserstein} of the Kantorovich-Rubinstein distance, it comes by choosing the trivial coupling $\mathcal{L}(\mu^N)\otimes\delta_u$ that
\begin{align}
\label{eq:W1muNdeltaU}
\mathcal{W}_1^{(\mathcal{P})}(\mathcal{L}(\mu^N), \delta_u)&\leq  \EE\Big[ 1\wedge 
\sup_{t\in[0,T]} \| \mu_{t}^N - u_{t} \|_{0} \Big] \nonumber\\
 &\leq \EE\Big[1\wedge \sup_{t\in[0,T]} \| \mu_{t}^N - u^N_{t}\|_{0} \Big] +   \EE\Big[1\wedge \sup_{t\in[0,T]} \sup_{ \|\phi\|_{L^\infty}\leq 1} \langle u_{t}^N - u_{t},\phi \rangle \Big] \nonumber \\
&\leq  \EE\Big[1\wedge \sup_{t\in[0,T]} \| \mu_{t}^N - u^N_{t}\|_{0} \Big] +  \EE\Big[1\wedge \sup_{t\in[0,T]}  \| u_{t}^N - u_{t}\|_{L^q(\T^d)} \Big].
\end{align}
We first treat the second term on the r.h.s. and will prove that 
\begin{equation}\label{eq:W1-uN}
    \EE\Big[1\wedge \sup_{t\in[0,T]}  \| u_{t}^N - u_{t}\|_{L^q(\T^d)} \Big] \leq C N^{-(\varrho-\varepsilon)}.
\end{equation}
Consider $N_0\equiv N_0(\omega)$ the smallest integer such that $\tau_{N_0}>T + \frac{1}{C_T} \log(2)$.
 Then we get
\begin{align*}
    \EE\Big[ 1\wedge \sup_{t\in[0,T]} \|u^N_t - u_t\|_{L^q(\mathbb{T}^d)} \Big] \leq  \EE\Big[ \mathbbm{1}_{\{N\geq N_0\}} \sup_{t\in[0,T]} \|u^N_t - u_t\|_{L^q(\mathbb{T}^d)} \Big] + \PP(N_0\geq N) .
\end{align*}
The above choice of $N_0$ induces that for $N\geq N_0$, we have that $\frac{2A_N}{1+A_N} e^{C_T T}= 2e^{C_T(T-\tau_N)} \leq 1$. Hence, in view of \eqref{eq:Bihari}, which holds true for any $t\leq T$ on the event $\{N\geq N_0\}$, we obtain
\begin{align*}
    \EE\Big[ \mathbbm{1}_{\{N\geq N_0\}} \sup_{t\in[0,T]} \|u^N_t - u_t\|_{L^q(\mathbb{T}^d)} \Big] &\leq e^{\frac{1}{3} C_T T} \EE\left[\Big( \frac{2A_N/(1+A_N)}{2 - e^{C_T T} 2A_N/(1+A_N )} \Big)^{\frac{1}{3}} \right]\\
    &\leq e^{\frac{1}{3} C_T T} \EE\left[\Big( \frac{2A_N}{1+A_N} \Big)^{\frac{1}{3}} \right]\\
    &\leq C \, e^{\frac{1}{3} C_T T}\, N^{-(\varrho-\varepsilon)},
\end{align*}
using the definition \eqref{eq:defAN} of $A_{N}$ and the fact that $A$ has finite moments.

\smallskip

Now we estimate $\PP(N_0\geq N)$. By the definition of $N_0$, we have that $\tau_{N_0-1} \leq T + \frac{1}{C_T}\log(2)$. Hence in view of the definition of $\tau_{N_0-1}$ and $A_N$, we deduce that
\begin{equation*}
    N_0 \leq 1 + \left( 2 e^{C_T T}-1\right)^{\frac{1}{3(\varrho-\varepsilon)}} A^{\frac{1}{\varrho-\varepsilon}} .
\end{equation*}
Now we get 
\begin{align*}
    \PP \left( N_0\geq N \right) &\leq \PP\left( A \geq \frac{(N-1)^{\varrho-\varepsilon}}{(2e^{C_T T}-1)^{\frac{1}{3}}} \right) \\
    & \leq C \frac{\EE A^p}{(N-1)^{p(\varrho-\varepsilon})} ,
\end{align*}
by the Markov inequality, for any $p\geq 1$. Hence \eqref{eq:W1-uN} follows.

\smallskip

For the first term in the right-hand side of \eqref{eq:W1muNdeltaU}, we observe that
\begin{align*}
|\langle \mu_{t}^{N}, \phi\rangle - \langle u^N_t, \phi \rangle| &=| \langle \mu_{t}^{N}, (\phi-\phi\ast V^{N})\rangle |\\
&\leq \Big\langle \mu_{t}^{N}, \int_{\T^{d}} V(y)~ |\phi(.)-  \phi( \frac{y}{N^{\alpha}}-.) | \,  dy \Big\rangle \\
&\leq \frac{C \|\phi\|_{\text{Lip}}}{N^{\alpha}}.
\end{align*}
In view of \eqref{eq:W1muNdeltaU}, \eqref{eq:W1-uN} the above inequality, and since $\alpha\geq \varrho$, the desired result follows.
\end{proof}

\begin{remark}\label{rk:BesselforKS}
Like for Burgers equation, it is possible to have convergence in Bessel spaces for the Keller-Segel model. We outline briefly how the above proof can be adapted to get a rate of convergence in $H^\beta_{p}(\T^d)$, assuming $p>d$, $\beta \in (\frac{d}{p},1)$ and $u_{0}\in H^\beta_{p}(\T^d)$. The latter point ensures that $u$ is a mild solution of the Keller-Segel PDE in $\mathcal{C}([0,T];H^\beta_{p}(\T^d))$, for some $T>0$.
\begin{itemize}
\item First, by a simple adaptation of the proof of Lemma~\ref{lem:boundHolder}, using the Calder\'on-Zygmund property of the Keller-Segel kernel, one can prove that for $f\in H^\beta_{p}(\T^d)$, $\lVert K\ast f \rVert_{\beta+1,p} \lesssim \lVert f \rVert_{\beta,p}$. By Sobolev embedding, this implies that in this case, $K\ast f$ is Lipschitz continuous and bounded.

\item Then mimicking the decomposition of $u-u^N$ as in the proof of Theorem~\ref{th:burgers}, one obtains exactly \eqref{eq:decompuN-u_Burgers} for the Keller-Segel equation, using in particular that $K\in L^1(\T^d)$ and therefore that $ \lVert K\ast f\rVert_{\beta,p} \leq \lVert K\rVert_{L^1(\T^d)} \lVert f\rVert_{\beta,p}$.

\item The commutator term $E$ is still given by \eqref{eq:defE-KS} but now for the $H^\beta_{p}$ norm, using the Lipschitz property established in the first point and the same proof as in the Burgers case, it reads
\begin{align*}
\lVert E_{t} \rVert_{\beta,p} &\lesssim \int_{0}^t \frac{1}{(t-s)^{\frac{1+\beta}{2}}} \lVert \langle \mu^N_{s},  V^N (x-\cdot) \left| K\ast u^N_s(\cdot)-  K\ast u^N_s(x)\right|\rangle \rVert_{L^p(\T^d)} \, ds\\
&\lesssim \frac{1}{N^\alpha} \int_{0}^t \frac{1}{(t-s)^{\frac{1+\beta}{2}}} \lVert K\ast u^N_{s}\rVert_{\text{Lip}} \lVert u^N_{s} \rVert_{L^p(\T^d)} \, ds\\
 &\lesssim \frac{1}{N^\alpha} \Big( 1 + \int_{0}^t \frac{1}{(t-s)^{\frac{1+\beta}{2} } } \|u^N_{s}-u_{s}\|_{\beta, p}^2\, ds \Big).
\end{align*}

\item As for the estimate on the stochastic integral, we have as in \eqref{eq:boundMB2} that $\big\|\sup_{s\in [0,t]} \| M^N_{s}\|_{\beta,p} \big\|_{L^m(\Omega)} \lesssim N^{-\frac{1}{2}(1-\alpha(2d+ 2\beta-2d/p) + \varepsilon}$.

\end{itemize}
Eventually this leads to a result similar to Theorem~\ref{TMain} with a rate $\tilde{\varrho}= \min \left(\alpha , \frac{1}{2}- \alpha (\beta+d(1-\frac{1}{p}))\right)$.
\end{remark}

\subsection{Trajectorial propagation of chaos for Keller-Segel}\label{sec:PoC-KS}

As in Section~\ref{subsec:trajectorialPoC-Burgers}, we tackle here the question of trajectorial convergence of the particles towards typical independent particles, which are solutions to the McKean-Vlasov equations:
\begin{align}\label{eq:McKVKS}
d\overline{X}_{t}^i = \frac{1}{2} K\ast u_{t}(\overline{X}_{t}^i)\, dt + dW^i_{t}, \quad \mathcal{L}(\overline{X}_{t}^i) = u_{t},
\end{align}
with the same Brownian motions and initial condition $X_{0}^i$ as the particle system. The above equation has a unique strong solution for the same reasons as \eqref{eq:McKVBurgers} noticing that uniqueness of the solution to the Keller-Segel equation holds (see Theorem~\ref{prop:PDE}) and the fact that  $K\ast u \in L^\infty([0,T];L^\infty(\T^d))$. We thus have the following result on a  quantitative convergence of our particle system towards \eqref{eq:McKVKS}:
\begin{proposition}\label{prop:KSPoC-particles}
With the previous notations and assuming that \eqref{H} holds true, then for any $\varepsilon>0$,
\begin{equation*}
\lim_{N\to +\infty} N^{\varrho-\varepsilon}  \sup_{t\in [0,T]} |X_t^{i,N}-\overline{X}_t^{i}| = 0 \ a.s.
\end{equation*}
\end{proposition}
To prove the above, it suffices to adapt line by line Proposition~\ref{prop:BurgersPoC-particles}, using Proposition~\ref{prop:Kelelr-Seguel} instead of Proposition~\ref{prop:BurgersLipschitz}. As before, the previous proposition  implies 
 the convergence of the empirical measure on the space of trajectories.
\begin{corollary}\label{cor:trajectorialPoC-Burgers}
The sequence $(\mu^N)_{N\in \N^*}$ converges in law in $\mathcal{P}(\mathcal{C}([0,T];\T^d))$ towards $u$.
\end{corollary}

\begin{remark}
As in Remark~\ref{rk:Burgers}, one could be interested in a notion of convergence in Wasserstein distance that is stronger than the one given in Theorem~\ref{th:burgers}(ii) and Theorem~\ref{TMain}(ii), namely a convergence at the level of the (random) empirical measure rather than its law. Indeed, at least formally, there is $\mathcal{W}^{(\mathcal{P})}_{1}(\mathcal{L}(\mu^N),\delta_{u}) 
\leq \EE[ \mathcal{W}_{1}(\mu^N ,u)]$. In Remark~\ref{rk:Burgers}, we were able to give a Wasserstein bound on the marginals. We could follow the same path here. We also notice that one can proceed as in \cite[Corollary 2.3]{ACM} to get the following bound in Kantorovich-Rubinstein distance:
\begin{align*}
\forall N\in \N^*,\quad \Big(\EE\Big[\sup_{t\in [0,T]} \|\mu^N_{t}-u_{t}\|_{0}^m\Big]\Big)^{1/m} \leq C\, N^{-\varrho+\varepsilon}.
\end{align*}
This bound requires to have convergence in $L^1$-norm in space of $\mu^N_{t}-u_{t}$, so that it can be derived here, but not in the Burgers case. However, as in Remark~\ref{rk:Burgers}, this bound does not imply the bound on $\mathcal{W}^{(\mathcal{P})}_{1}(\mathcal{L}(\mu^N),\delta_{u})$ from Theorem~\ref{TMain}(ii).
\end{remark}

\paragraph{Acknowledgments.} 
This work is supported by the SDAIM project ANR-22-CE40-0015 jointly funded by the S\~ao Paulo Research Foundation (FAPESP) and the French National Research Agency (ANR).\\
The authors thank two anonymous referees for their comments on the first version of the paper which have led us to improve the overall presentation of the paper, detail the trajectorial propagation of chaos result and find suitable quantitative estimates for the initial data.

\appendix

\section{Stochastic convolution integrals}\label{App}

The goal of this appendix is to prove the estimates \eqref{eq:boundMB2} and \eqref{eq:boundMN-KS} on the stochastic integrals defined respectively in \eqref{eq:defMB} for Burgers and in \eqref{eq:defMN-KS} for Keller-Segel. Since we will deal here with both $\mathcal{D}=\R$ and $\mathcal{D}=\T^d$, we denote by $\|\cdot\|_{H^\beta_{p}(\mathcal{D})}$ the Bessel norm in $H^\beta_{p}(\mathcal{D})$, which was defined in \eqref{eq:Besselnorm} on $\R$ and in \eqref{eq:Besselnorm-torus} on the torus.

In Proposition B.3 of \cite{ACM}, a bound on the $H^\beta_{2}(\R^d)$ norm of $M^N$ was proven. The differences compared to this work are that: 
\begin{enumerate}[label=(\roman*)]
\item here the particles are defined without cut-off.

\item In addition, in the Burgers case the kernel is not covered by the assumptions of \cite{ACM};

\item in the Keller-Segel case we now work on the torus.
\end{enumerate}   
In the proof of \cite[Proposition B.3]{ACM}, the nature of the particles is unimportant as they play no role in the computations. So points (i) and (ii) will not affect the result. 
As for the choice of the space, since the proof of \cite[Proposition B.3]{ACM} relies essentially on the Burkholder-Davis-Gundy inequality on Hilbert spaces (which $H^\beta_{2}(\T^d)$ still is) and Sobolev embeddings that are still valid on the torus, one can check that all computations go through similarly on $\T^d$ instead of $\R^d$.
Hence we can state the following bound on $M^N$, which is similar to the one from \cite[Proposition B.3]{ACM}.
\begin{proposition}
\label{prop:martingale-bound-Hbeta}
Let $\beta \in \R$, $m\geq 1$ and $\alpha\in [0,1]$ that defines $V^N$ as in \eqref{eq:defVN}. For any $\delta\in(0,1]$, there exists $C>0$ such that
\begin{align*}
\left\| \|M^N_{t}-M^N_{s}\|_{H^\beta_{2}(\mathcal{D})} \right\|_{L^m(\Omega)}  \leq C\, (t-s)^{\frac{\delta}{2}} \, N^{-\frac{1}{2} \left( 1-\alpha(d + 4\delta+2\beta) \right)} ,\quad \forall s\leq t\in[0,T],~\forall N\in\N^*.
\end{align*}
\end{proposition}

We now give a version of Garsia-Rodemich-Rumsey's Lemma \cite{GRR} (for $\R$-valued processes, the following lemma already appears in \cite[Corollary 4.4]{RTY}, and the extension to Banach spaces is consistent with Garsia-Rodemich-Rumsey's Lemma with no additional difficulty, see e.g. \cite[Theorem A.1]{FrizVictoir}).
	\begin{lemma}\label{lem:GRR}
		Let $E$ be a Banach space and $(Y^n)_{n\ge 1}$ be a sequence of $E$-valued continuous processes on $[0,T]$.
		Let $m\geq 1$ and $\eta>0$ such that $m\eta>1$ and assume that there exists a constant $C_{0} > 0$
		and a sequence $(\delta_{n})_{n\ge 1}$ of positive real numbers such that
		\begin{equation*}
			\left(\EE \Big[ \big\|Y^n_{s} - Y^n_{t} \big\|_{E}^m \Big] \right)^\frac{1}{m}
			 \leq 
			C_{0}  |s-t|^\eta\, \delta_{n}, \quad \forall s,t\in[0,T], ~\forall n \ge 1 .
		\end{equation*}
		Then for any $m_{0}\in (0,m]$, there exists a constant 
		$C$,
		 depending only on $C_{0}$, $m$, $m_{0}$, $\eta$ and $T$, such that $\forall n \ge 1$,
		\begin{equation*}
			\left(\EE \Big[ \sup_{t\in[0,T]} \big\| Y^n_t-Y^n_0 \big\|_{E}^{m_{0}} \Big] \right)^\frac{1}{m_{0}}
			\leq
			C ~ \delta_{n}.
		\end{equation*}
	\end{lemma}

These two results can be combined to obtain the following bound:
\begin{proposition}
\label{prop:martingale-bound}
Let $\beta\in \R$, $m\geq 1$ and $\alpha\in [0,1]$ that defines $V^N$ as in \eqref{eq:defVN}.
 Let $\varepsilon \in (0,2\alpha)$. Then there exists $C>0$ such that for any $t\in [0,T]$ and $N\in\N^*$,
\begin{equation*}
\Big\| \sup_{s\in [0,t] } \|M^N_{s}\|_{H^\beta_{2}(\mathcal{D})} \Big\|_{L^{m}(\Omega)}  \leq C\,  N^{- \frac{1}{2}\left(1-\alpha (d+ 2\beta) \right) + \varepsilon} ,
\end{equation*}

\end{proposition}

\begin{proof}
	We aim to apply Lemma \ref{lem:GRR} to  $M^N$ in the Hilbert space $H^\beta_2(\D)$.
	Let  $\varepsilon >0$ and $m_{0}>0$. With the notations of Proposition \ref{prop:martingale-bound-Hbeta}, 
	let us choose $\delta = \frac{\varepsilon}{2\alpha}$,  
	 $\eta = \tfrac{\delta}{2}$ and $\delta_{N} = N^{-\rho}$ with $\rho = - \frac{1}{2}\left(1-\alpha (1+ 2\beta) \right)  + 2 \alpha\delta = - \frac{1}{2}\left(1-\alpha (1+ 2\beta) \right)  + \varepsilon$. Hence, 
	choosing $m \geq 1\vee m_{0}$ large enough so that $ m \eta >1$,
	the inequality in Proposition \ref{prop:martingale-bound-Hbeta} shows that $M^N$ satisfies the conditions of Lemma \ref{lem:GRR} and the 
	 desired result follows.
\end{proof}


\section{Quantitative approximation of a density in Bessel norm by an i.i.d. sample}\label{app:initialLLN}

The aim of this appendix is to prove the following result, which quantifies, in Bessel norm, the approximation of a density $u_{0}$ by the smoothed empirical measure of i.i.d. variables distributed according to $u_{0}$.

\begin{proposition}\label{prop:initialrate}
Let $\D$ be either $\T^d$ or $\R^d$. Let $\beta'\geq 0$, $p\in [2,+\infty)$ and $u_{0}\in H^{\beta'}_{p}(\D)$ be a probability density. 
Let $X_{i} \sim u_{0}$ be a family of i.i.d. random variables and denote $u_{0}^N(x) = \frac{1}{N} \sum_{k=1}^N V^N(x-X_{i})$ the smoothed empirical measure, where $V^N$ is the smooth approximation of identity from \eqref{eq:defVN}.  If $\D = \R^d$, assume that $u_{0}$ is compactly supported and that $V$ can be written as $V=W\ast W$, where $W$ is again a smooth compactly supported probability density function. 
Then for $\beta \leq \beta'$, $m\geq p$,
\begin{equation*}
\EE\left[  \left\|  u^N_{0} - u_{0} \right\|_{\beta,p}^m   \right]^{1/m} \lesssim N^{-\alpha(\beta'-\beta)} + N^{\alpha\beta + d\alpha(1-\frac{2}{m}) - \frac{1}{2}},\quad \forall N\geq 1.
\end{equation*}
\end{proposition}

\begin{proof}
The error is split in two as follows
\begin{equation*}
\EE\left[  \left\|  u^N_{0} - u_{0} \right\|_{\beta,p}^m   \right]^{1/m} \leq \EE\left[  \left\|  u^N_{0} - \EE u_{0}^N \right\|_{\beta,p}^m   \right]^{1/m} +  \left\|  \EE u^N_{0} - u_{0} \right\|_{\beta,p}.
\end{equation*}
As $\EE u^N_{0} = V^N \ast u_{0}$, a classical approximation estimate yields for the second term that 
\begin{equation}\label{eq:boundinit1}
\|V^N \ast u_{0} - u_{0} \|_{\beta,p} \lesssim N^{-\alpha(\beta'-\beta)} \lVert u_{0} \rVert_{\beta',p}.
\end{equation}

Let us now focus on the first term in the RHS of the previous inequality.
There is
\begin{align*}
\EE\left[  \left\|  u^N_{0} - \EE u_{0}^N \right\|_{\beta,p}^m   \right] 
&= \EE\Big[ \Big( \int_{\D} \big| \frac{1}{N} \sum_{i=1}^N (I-\Delta)^{\frac{\beta}{2}} \big\{V^N(\cdot-X_{i}) - V^N\ast u_{0}(\cdot) \}(x) \big|^p \, dx \Big)^\frac{m}{p} \Big].\end{align*}
In the case $\D=\T^d$ the compactness allows to proceed with Jensen's inequality without further assumption on $V$. In general, we instead rely on $V^N = W^N\ast W^N$, with $W^N(x) = N^{d\alpha} W(N^\alpha x)$, to write
\begin{align*}
\EE\left[  \left\|  u^N_{0} - \EE u_{0}^N \right\|_{\beta,p}^m \right] 
&= \EE\Big[ \Big( \int_{\D} \big| \frac{1}{N} \sum_{i=1}^N (I-\Delta)^{\frac{\beta}{2}}W \ast \big\{W^N(\cdot-X_{i}) - W^N\ast u_{0}(\cdot) \}(x) \big|^p \, dx \Big)^\frac{m}{p} \Big] \\
&\leq \lVert (I-\Delta)^{\frac{\beta}{2}}W^N \rVert_{L^1(\D)}^m\, \EE\Big[ \Big( \int_{\D} \big| \frac{1}{N} \sum_{i=1}^N \big\{ W^N(x-X_{i}) - W^N\ast u_{0}(x) \big\} \big|^p \, dx \Big)^\frac{m}{p} \Big] ,
\end{align*}
where we used a convolution inequality in the last line. 
Then, since $\lVert (I-\Delta)^{\frac{\beta}{2}}W^N \rVert_{L^1(\D)}^m \lesssim N^{\alpha \beta m}$ and by Jensen's inequality, we have
\begin{align}\label{eq:boundJensen}
\EE\left[  \left\|  u^N_{0} - \EE u_{0}^N \right\|_{\beta,p}^m   \right] 
&\leq C\, N^{\alpha \beta m}\, \EE \int_{\D} \big| \frac{1}{N} \sum_{i=1}^N \big\{ W^N(x-X_{i}) - W^N\ast u_{0}(x) \big\} \big|^m \, dx ,
\end{align}
where $C$ depends on $W$, $\beta$, $m$, $p$ and the supports of $W$ and $u_{0}$, but not on $N$.

Without loss of generality, assume that $m$ is an even integer (as if it is not, one can always choose $m$ larger to be an even integer). Let $k\in \N^*$ such that $m=2k$. Since $W^N\ast u_{0}(x) = \EE[W^N(x-X_{i})]$, we notice that $\sum_{i=1}^N \big\{ W^N(x-X_{i}) - W^N\ast u_{0}(x) \big\}$ is a centered random walk, thus we know that there is a constant $C$ such that for any $N$, 
\begin{align*}
\EE \Big| \sum_{i=1}^N \big\{ W^N(x-X_{i}) - W^N\ast u_{0}(x) \big\} \Big|^{2k} \leq C  \, N^k\, (2k-1)!!\ \EE  \big| W^N(x-X_{i}) - W^N\ast u_{0}(x) \big|^{2k}.
\end{align*}
Plugging the previous bound into \eqref{eq:boundJensen} for $m=2k$ yields
\begin{align*}
\EE\left[  \left\|  u^N_{0} - \EE u_{0}^N \right\|_{\beta,p}^m \right] 
&\lesssim N^{\alpha \beta m}\, N^{-\frac{m}{2}} \int_{\D}\int_{\D} \big| W^N(x-y) - W^N\ast u_{0}(x) \big|^{m} u_{0}(y)\, dy\, dx\\
&\lesssim N^{m(\alpha\beta-\frac{1}{2})} \int_{\D}\int_{\D} \Big| \int_{\D} \big(W^N(x-y) - W^N(x-z)\big) u_{0}(z)\, dz \Big|^{m} u_{0}(y)\, dy\, dx .
\end{align*}
Now Jensen's inequality and a change of variables give
\begin{align*}
\EE\left[  \left\|  u^N_{0} - \EE u_{0}^N \right\|_{\beta,p}^m \right] 
&\lesssim N^{m(\alpha\beta-\frac{1}{2})} \int_{\D^3} \big|W^N(x-y) - W^N(x-z)\big|^{m} u_{0}(z)\,  u_{0}(y) \, dz \, dy\, dx\\
&\lesssim N^{m(\alpha\beta-\frac{1}{2})} N^{d\alpha(m-2)} \int_{\D^3} |W(y) - W(z)|^m \, u_{0}(x+N^{-\alpha}z)\,  u_{0}(x+N^{-\alpha} y) \, dz \, dy\, dx .
\end{align*}
Having assumed that $p\geq 2$, we have by Cauchy-Schwarz inequality that $\int_{\D} u_{0}(x+N^{-\alpha}z)\,  u_{0}(x+N^{-\alpha} y) \, dx \leq \lVert u_{0}\rVert_{L^2(\D)}^2$, so that finally
\begin{equation}\label{eq:boundinit2}
\left(\EE\left[  \left\|  u^N_{0} - \EE u_{0}^N \right\|_{\beta,p}^m \right] \right)^{\frac{1}{m}} \lesssim N^{\alpha\beta-\frac{1}{2} + d\alpha(1-\frac{2}{m})} .
\end{equation}
Putting together \eqref{eq:boundinit1} and \eqref{eq:boundinit2} gives the result.
\end{proof}

 \section{Improved space regularity of the solutions of the PDEs}\label{app:regPDE}
 
 In this appendix, we study the space regularity of the Burgers and Keller-Segel equations. We obtain improved space regularity compared respectively to Theorem~\ref{th:PDEB} and Proposition~\ref{prop:PDE}, at the price of a singularity in time near $0$.
 
 First, consider the regularity for the Burgers equation.
 \begin{proposition}\label{prop:BurgersLipschitz}
Let $p > 1$, $\beta \in (1/p,1)$, $T>0$, $u_0 \in  L^1\cap H^{\beta}_{p}(\R)$, and let $u$ be the mild solution to the Burgers equation~\eqref{Burgers} given by Theorem~\ref{th:PDEB}. For any $\delta \in (0, 1+\beta)$,  there exist $C_{T} >0$ such that 
 \[
\forall t\in (0,T],\quad \| u_{t} \|_{ 1+ \delta,p} \leq \frac{ C_{T}}{t^{(1+\delta-\beta)/2}}. 
 \]
\end{proposition}
  
\begin{proof} 
 As the unique mild solution to \eqref{Burgers}, $u$ verifies 
 \[
 u_{t} =  e^{t\Delta} u_0 - \frac{1}{2} \int_0^t \partial_{x} e^{(t-s)\Delta }  u_{s}^{2}\ ds, \quad 0 \leq t \leq T.
\] 
Applying the Bessel norm and the triangular inequality we get 
 \[
\| u_{t}\|_{1+ \delta, p} 
 \leq \| e^{t\Delta} u_0\|_{1+ \delta, p}   + \frac{1}{2} \int_0^t \|\partial_{x} e^{(t-s)\Delta }  u_{s}^{2} \|_{1+ \delta, p} \ ds.
\] 
Observe that
\begin{equation}\label{initial}
\| e^{t\Delta} u_0\|_{1+ \delta, p} \leq  \frac{C}{t^{(1 + \delta - \beta)/2}},  
\end{equation}
 and 
\begin{align}\label{nonli}
  \frac{1}{2} \int_0^t \|\partial_{x} e^{(t-s)\Delta }  u_{s}^{2} \|_{1+ \delta, p} \ ds
&\leq  C \int_0^t  \frac{1}{(t-s)^{(1+\delta-\beta)/2}} \ \| u_{s}^{2} \|_{\beta, p}  \  ds \nonumber \\
& \leq  \frac{1}{2} \int_0^t    \frac{1}{(t-s)^{(1+\delta-\beta)/2}} \  \| u_{s} \|_{\beta, p}^{2}     \ ds\leq C_{T}  ,
\end{align}
 where we used convolution inequality, Lemma \ref{lemma:heatBeta}, Lemma  \ref{lem:equivNorms} and that 
 $H^{\beta}_{p}(\mathbb{R}^{d})$ is an algebra. 
We deduce the result from \eqref{initial} and \eqref{nonli}.
\end{proof}

 An adaptation of the previous proof yields the following result for the Keller-Segel PDE. 
 \begin{proposition}\label{prop:Kelelr-Seguel}
Let $q > d$ , $T>0$, $u_0 \in  L^1\cap L^{q}(\T^{d})$, and let $u$ be the mild solution to the Keller equation  \eqref{eq:mildKS} given by Proposition 
\ref{prop:PDE}. For any $\delta \in (0,1)$,  there exist $C_{T} >0$ such that 
 \[
\forall t\in (0,T],\quad \| u_{t} \|_{ 1+ \delta,q} \leq \frac{ C_{T}}{t^{(1+\delta)/2}}. 
 \]
\end{proposition}


\begin{thebibliography}{31}
\providecommand{\natexlab}[1]{#1}
\providecommand{\url}[1]{\texttt{#1}}
\expandafter\ifx\csname urlstyle\endcsname\relax
  \providecommand{\doi}[1]{doi: #1}\else
  \providecommand{\doi}{doi: \begingroup \urlstyle{rm}\Url}\fi

\bibitem[Bahouri et~al.(2011)Bahouri, Chemin, and Danchin]{BCD}
H.~Bahouri, J.-Y. Chemin, and R.~Danchin.
\newblock \emph{Fourier analysis and nonlinear partial differential equations},
  volume 343 of \emph{Grundlehren der mathematischen Wissenschaften
  [Fundamental Principles of Mathematical Sciences]}.
\newblock Springer, Heidelberg, 2011.

\bibitem[Biler(2018)]{Biler}
P.~Biler.
\newblock Mathematical challenges in the theory of chemotaxis.
\newblock \emph{Ann. Math. Sil.}, 32\penalty0 (1):\penalty0 43--63, 2018.

\bibitem[Bogachev(2007)]{BogachevII}
V.~I. Bogachev.
\newblock \emph{Measure theory. {V}ol. {II}}.
\newblock Springer-Verlag, Berlin, 2007.

\bibitem[Bossy and Talay(1996)]{BossyTalay0}
M.~Bossy and D.~Talay.
\newblock Convergence rate for the approximation of the limit law of weakly
  interacting particles: application to the {B}urgers equation.
\newblock \emph{Ann. Appl. Probab.}, 6\penalty0 (3):\penalty0 818--861, 1996.

\bibitem[Bossy and Talay(1997)]{BossyTalay}
M.~Bossy and D.~Talay.
\newblock A stochastic particle method for the {M}c{K}ean-{V}lasov and the
  {B}urgers equation.
\newblock \emph{Math. Comp.}, 66\penalty0 (217):\penalty0 157--192, 1997.

\bibitem[Bresch et~al.(2023)Bresch, Jabin, and Wang]{BJW2020}
D.~Bresch, P.-E. Jabin, and Z.~Wang.
\newblock Mean field limit and quantitative estimates with singular attractive
  kernels.
\newblock \emph{Duke Math. J.}, 172\penalty0 (13):\penalty0 2591--2641, 2023.

\bibitem[Calderoni and Pulvirenti(1983)]{Calde}
P.~Calderoni and M.~Pulvirenti.
\newblock Propagation of chaos for {B}urgers' equation.
\newblock \emph{Ann. Inst. H. Poincar\'{e} Sect. A (N.S.)}, 39\penalty0
  (1):\penalty0 85--97, 1983.

\bibitem[Cattiaux and P{\'e}d\`eches(2016)]{CattiauxPedeches}
P.~Cattiaux and L.~P{\'e}d\`eches.
\newblock {The 2-D stochastic Keller-Segel particle model: existence and
  uniqueness}.
\newblock \emph{ALEA}, 13:\penalty0 447--463, 2016.

\bibitem[Cazacu(2025)]{Cazacu}
N.~Cazacu.
\newblock {Stochastic numerical approximation for nonlinear Fokker-Planck
  equations with singular kernels}.
\newblock \emph{Preprint  \href{https://arxiv.org/abs/2504.06132}{arXiv:2504.06132}} , 2025.

\bibitem[Chen et~al.(2022)Chen, Holzinger, J\"{u}ngel, and
  Zamponi]{HolzingerEtAl-Porous}
L.~Chen, A.~Holzinger, A.~J\"{u}ngel, and N.~Zamponi.
\newblock Analysis and mean-field derivation of a porous-medium equation with
  fractional diffusion.
\newblock \emph{Comm. Partial Differential Equations}, 47\penalty0
  (11):\penalty0 2217--2269, 2022.

\bibitem[Chen et~al.(2024)Chen, Holzinger, and J{\"u}ngel]{HolzingerEtAl}
L.~Chen, A.~Holzinger, and A.~J{\"u}ngel.
\newblock {Fluctuations around the mean-field limit for attractive Riesz
  potentials in the moderate regime}.
\newblock \emph{Preprint
  \href{https://arxiv.org/abs/2405.15128}{arXiv:2405.15128}}, 2024.

\bibitem[Chodron~de Courcel et~al.(2023)Chodron~de Courcel, Rosenzweig, and
  Serfaty]{ChodronEtAl}
A.~Chodron~de Courcel, M.~Rosenzweig, and S.~Serfaty.
\newblock {Sharp uniform-in-time mean-field convergence for singular periodic
  Riesz flows}.
\newblock \emph{Ann. Inst. H. Poincar\'e C Anal. Non Lin\'eaire (to appear),
  \href{https://arxiv.org/abs/2304.05315}{arXiv:2304.05315}}, 2023.

\bibitem[Cole(1951)]{Cole}
J.~D. Cole.
\newblock On a quasi-linear parabolic equation occurring in aerodynamics.
\newblock \emph{Quart. Appl. Math.}, 9:\penalty0 225--236, 1951.

\bibitem[Dragomir(2003)]{Dragomir}
S.~S. Dragomir.
\newblock \emph{Some {G}ronwall type inequalities and applications}.
\newblock Nova Science Publishers, Inc., Hauppauge, NY, 2003.

\bibitem[Flandoli and Leocata(2019)]{FlandoliLeocata}
F.~Flandoli and M.~Leocata.
\newblock A particle system approach to aggregation phenomena.
\newblock \emph{J. Appl. Probab.}, 56\penalty0 (1):\penalty0 282--306, 2019.

\bibitem[Flandoli et~al.(2019)Flandoli, Leimbach, and
  Olivera]{FlandoliLeimbachOlivera}
F.~Flandoli, M.~Leimbach, and C.~Olivera.
\newblock {Uniform convergence of proliferating particles to the FKPP
  equation}.
\newblock \emph{J. Math. Anal. Appl.}, 473\penalty0 (1):\penalty0 27--52, 2019.

\bibitem[Flandoli et~al.(2020)Flandoli, Olivera, and
  Simon]{FlandoliOliveraSimon}
F.~Flandoli, C.~Olivera, and M.~Simon.
\newblock Uniform approximation of 2 dimensional {N}avier-{S}tokes equation by
  stochastic interacting particle systems.
\newblock \emph{SIAM J. Math. Anal.}, 52\penalty0 (6):\penalty0 5339--5362,
  2020.

\bibitem[Fournier and Guillin(2015)]{FournierGuillin}
N.~Fournier and A.~Guillin.
\newblock On the rate of convergence in {W}asserstein distance of the empirical
  measure.
\newblock \emph{Probab. Theory Related Fields}, 162\penalty0 (3-4):\penalty0
  707--738, 2015.

\bibitem[Fournier and Jourdain(2017)]{FournierJourdain}
N.~Fournier and B.~Jourdain.
\newblock {Stochastic particle approximation of the Keller--Segel equation and
  two-dimensional generalization of Bessel processes}.
\newblock \emph{Ann. Appl. Probab.}, 27\penalty0 (5):\penalty0 2807--2861,
  2017.

\bibitem[Fournier and Tardy(2024)]{FournierTardy}
N.~Fournier and Y.~Tardy.
\newblock Collisions of the supercritical {K}eller-{S}egel particle system.
\newblock \emph{J. Eur. Math. Soc. (to appear),
  \href{https://arxiv.org/abs/2110.08490}{arXiv:2110.08490}}, 2024.

\bibitem[Friz and Victoir(2010)]{FrizVictoir}
P.~K. Friz and N.~B. Victoir.
\newblock \emph{Multidimensional stochastic processes as rough paths}, volume
  120 of \emph{Cambridge Studies in Advanced Mathematics}.
\newblock Cambridge University Press, Cambridge, 2010.

\bibitem[Garsia et~al.(1970/71)Garsia, Rodemich, and Rumsey]{GRR}
A.~M. Garsia, E.~Rodemich, and H.~Rumsey, Jr.
\newblock A real variable lemma and the continuity of paths of some {G}aussian
  processes.
\newblock \emph{Indiana Univ. Math. J.}, 20:\penalty0 565--578, 1970/71.

\bibitem[Gilbarg and Trudinger(2001)]{GilbargTrudinger}
D.~Gilbarg and N.~S. Trudinger.
\newblock \emph{Elliptic partial differential equations of second order}.
\newblock Classics in Mathematics. Springer-Verlag, Berlin, 2001.
\newblock Reprint of the 1998 edition.

\bibitem[Gutkin and Kac(1983)]{Kac2}
E.~Gutkin and M.~Kac.
\newblock Propagation of chaos and the {B}urgers equation.
\newblock \emph{SIAM J. Appl. Math.}, 43\penalty0 (4):\penalty0 971--980, 1983.

\bibitem[Hao et~al.(2024)Hao, Jabir, Menozzi, R{\"o}ckner, and
  Zhang]{JabirEtAl}
Z.~Hao, J.-F. Jabir, S.~Menozzi, M.~R{\"o}ckner, and X.~Zhang.
\newblock {Propagation of chaos for moderately interacting particle systems
  related to singular kinetic Mckean-Vlasov SDEs}.
\newblock \emph{Preprint
  \href{https://arxiv.org/abs/2405.09195}{arXiv:2405.09195}}, 2024.

\bibitem[Jabin and Wang(2018)]{JabinWang}
P.-E. Jabin and Z.~Wang.
\newblock Quantitative estimates of propagation of chaos for stochastic systems
  with {$W^{-1,\infty}$} kernels.
\newblock \emph{Invent. Math.}, 214\penalty0 (1):\penalty0 523--591, 2018.

\bibitem[Jourdain and M\'{e}l\'{e}ard(1998)]{JourdainMeleard}
B.~Jourdain and S.~M\'{e}l\'{e}ard.
\newblock Propagation of chaos and fluctuations for a moderate model with
  smooth initial data.
\newblock \emph{Ann. Inst. H. Poincar\'{e} Probab. Statist.}, 34\penalty0
  (6):\penalty0 727--766, 1998.

\bibitem[Lemari\'{e}-Rieusset(2002)]{Lemai}
P.~G. Lemari\'{e}-Rieusset.
\newblock \emph{Recent developments in the {N}avier-{S}tokes problem}, volume
  431 of \emph{Chapman \& Hall/CRC Research Notes in Mathematics}.
\newblock Chapman \& Hall/CRC, Boca Raton, FL, 2002.

\bibitem[McKean(1967)]{McKean}
H.~P. McKean, Jr.
\newblock Propagation of chaos for a class of non-linear parabolic equations.
\newblock In \emph{Stochastic {D}ifferential {E}quations ({L}ecture {S}eries in
  {D}ifferential {E}quations, {S}ession 7, {C}atholic {U}niv., 1967)}, Lecture
  Series in Differential Equations, Session 7, pages 41--57. Air Force Office
  of Scientific Research, Office of Aerospace Research, United States Air
  Force, Arlington, VA, 1967.

\bibitem[M\'{e}l\'{e}ard and Roelly-Coppoletta(1987)]{Meleard}
S.~M\'{e}l\'{e}ard and S.~Roelly-Coppoletta.
\newblock {A propagation of chaos result for a system of particles with
  moderate interaction}.
\newblock \emph{Stochastic Process. Appl.}, 26\penalty0 (2):\penalty0 317--332,
  1987.

\bibitem[Oelschl\"{a}ger(1985)]{Oelschlager85}
K.~Oelschl\"{a}ger.
\newblock A law of large numbers for moderately interacting diffusion
  processes.
\newblock \emph{Z. Wahrsch. Verw. Gebiete}, 69\penalty0 (2):\penalty0 279--322,
  1985.

\bibitem[Oelschl\"{a}ger(1987)]{Oelschlager87}
K.~Oelschl\"{a}ger.
\newblock A fluctuation theorem for moderately interacting diffusion processes.
\newblock \emph{Probab. Theory Related Fields}, 74\penalty0 (4):\penalty0
  591--616, 1987.

\bibitem[Oelschl\"{a}ger(1990)]{OelschlagerPorous}
K.~Oelschl\"{a}ger.
\newblock Large systems of interacting particles and the porous medium
  equation.
\newblock \emph{J. Differential Equations}, 88\penalty0 (2):\penalty0 294--346,
  1990.

\bibitem[Olivera et~al.(2023)Olivera, Richard, and Toma\v{s}evi\'{c}]{ACM}
C.~Olivera, A.~Richard, and M.~Toma\v{s}evi\'{c}.
\newblock Quantitative particle approximation of nonlinear {F}okker-{P}lanck
  equations with singular kernel.
\newblock \emph{Ann. Sc. Norm. Super. Pisa Cl. Sci. (5)}, 24\penalty0
  (2):\penalty0 691--749, 2023.

\bibitem[Osada and Kotani(1985)]{Osada}
H.~Osada and S.~Kotani.
\newblock Propagation of chaos for the {B}urgers equation.
\newblock \emph{J. Math. Soc. Japan}, 37\penalty0 (2):\penalty0 275--294, 1985.

\bibitem[Perthame(2004)]{perthame}
B.~Perthame.
\newblock P{DE} models for chemotactic movements: parabolic, hyperbolic and
  kinetic.
\newblock \emph{Appl. Math.}, 49\penalty0 (6):\penalty0 539--564, 2004.

\bibitem[Richard et~al.(2021)Richard, Tan, and Yang]{RTY}
A.~Richard, X.~Tan, and F.~Yang.
\newblock Discrete-time simulation of stochastic {V}olterra equations.
\newblock \emph{Stochastic Process. Appl.}, 141:\penalty0 109--138, 2021.

\bibitem[Schmeisser and Triebel(1987)]{SchmeisserTriebel}
H.-J. Schmeisser and H.~Triebel.
\newblock \emph{Topics in {F}ourier analysis and function spaces}.
\newblock A Wiley-Interscience Publication. John Wiley \& Sons, Ltd.,
  Chichester, 1987.

\bibitem[Senba and Suzuki(2002)]{SenbaSuzuki}
T.~Senba and T.~Suzuki.
\newblock Weak solutions to a parabolic-elliptic system of chemotaxis.
\newblock \emph{J. Funct. Anal.}, 191\penalty0 (1):\penalty0 17--51, 2002.

\bibitem[Serfaty(2020)]{Serfaty}
S.~Serfaty.
\newblock Mean field limit for {C}oulomb-type flows.
\newblock \emph{Duke Math. J.}, 169\penalty0 (15):\penalty0 2887--2935, 2020.

\bibitem[Simon and Olivera(2018)]{Simon}
M.~Simon and C.~Olivera.
\newblock {Non-local conservation law from stochastic particle systems}.
\newblock \emph{J. Dynam. Differential Equations}, 30\penalty0 (4):\penalty0
  1661--1682, 2018.

\bibitem[Sznitman(1986)]{SzTi}
A.-S. Sznitman.
\newblock A propagation of chaos result for {B}urgers' equation.
\newblock \emph{Probab. Theory Relat. Fields}, 71\penalty0 (4):\penalty0
  581--613, 1986.

\bibitem[Tardy(2024)]{Tardy}
Y.~Tardy.
\newblock Weak convergence of the empirical measure for the {K}eller-{S}egel
  model in both subcritical and critical cases.
\newblock \emph{Electron. J. Probab.}, 29:\penalty0 Paper No. 142, 35, 2024.

\bibitem[Triebel(1978)]{Triebel}
H.~Triebel.
\newblock \emph{Interpolation theory, function spaces, differential operators},
  volume~18 of \emph{North-Holland Mathematical Library}.
\newblock North-Holland Publishing Co., Amsterdam-New York, 1978.

\bibitem[Triebel(1983)]{TriebelTh}
H.~Triebel.
\newblock \emph{Theory of function spaces}, volume~78 of \emph{Monographs in
  Mathematics}.
\newblock Birkh\"{a}user Verlag, Basel, 1983.

\bibitem[van Neerven et~al.(2007)van Neerven, Veraar, and Weis]{vNeervenEtAl}
J.~M. A.~M. van Neerven, M.~C. Veraar, and L.~Weis.
\newblock Stochastic integration in {UMD} {B}anach spaces.
\newblock \emph{Ann. Probab.}, 35\penalty0 (4):\penalty0 1438--1478, 2007.

\bibitem[Veretennikov(1980)]{Veretennikov}
A.~J. Veretennikov.
\newblock Strong solutions and explicit formulas for solutions of stochastic
  integral equations.
\newblock \emph{Mat. Sb. (N.S.)}, 111(153)\penalty0 (3):\penalty0 434--452,
  480, 1980.

\bibitem[Wang et~al.(2023)Wang, Zhao, and Zhu]{WangEtAl}
Z.~Wang, X.~Zhao, and R.~Zhu.
\newblock Gaussian fluctuations for interacting particle systems with singular
  kernels.
\newblock \emph{Arch. Ration. Mech. Anal.}, 247\penalty0 (5):\penalty0 Paper
  No. 101, 62, 2023.

\end{thebibliography}
\end{document}